\documentclass[12pt,english]{article}

\usepackage[latin9]{inputenc}
\usepackage{amsthm}
\usepackage{amsmath}
\usepackage{amssymb}
\usepackage{graphicx}
\usepackage{setspace}
\usepackage{esint}

\PassOptionsToPackage{normalem}{ulem}
\usepackage{ulem}
\usepackage{tabularx}
\usepackage{mathrsfs}
\usepackage{subfigure}
\usepackage{longtable}
\usepackage{lscape}
\usepackage{bm}
\usepackage{color}
\usepackage{filecontents}

\makeatletter

\providecommand{\tabularnewline}{\\}

\numberwithin{equation}{section}
\numberwithin{figure}{section}
  \theoremstyle{definition}
  \newtheorem{condition}{\protect\conditionname}
\theoremstyle{plain}
\newtheorem{thm}{\protect\theoremname}
  \theoremstyle{remark}
  \newtheorem{rem}{\protect\remarkname}
  \theoremstyle{plain}
  \newtheorem{prop}{\protect\propositionname}
 \theoremstyle{definition}
 \newtheorem*{defn*}{\protect\definitionname}

\usepackage{graphics}\usepackage{amsfonts}\usepackage{amstext}\usepackage{tabularx}\usepackage{mathrsfs}\usepackage{subfigure}\usepackage{lscape}\usepackage{bm}

\newcommand{\sfrac}[2]{\ensuremath\raisebox{1.5pt}{\footnotesize
             $#1$}\kern-1pt/\kern-1pt 
             \raisebox{-2pt}{\footnotesize $#2$}}

\makeatother

\usepackage{babel}
  \providecommand{\conditionname}{Condition}
  \providecommand{\definitionname}{Definition}
  \providecommand{\propositionname}{Proposition}
  \providecommand{\remarkname}{Remark}
\providecommand{\theoremname}{Theorem}

\usepackage{natbib}

\begin{document}

\title{An Empirical Likelihood-based Local Estimation%
\thanks{The author would like to express his appreciation to Peter Boswijk,
Kees Jan van Garderen, Yuichi Kitamura, Richard Smith, Kenneth Judd,
Paulo Parente, participants in the seminars at University of Warwick,
Toulouse School of Economics, Cowles Foundation, University of Amsterdam
and participants in Econometric Society World Congress at Shanghai,
Winter North American Econometric Society at Atlanta for helpful comments
and discussions. All the remaining errors are mine.%
}}

\author{Zhengyuan Gao%
\thanks{The University of Iowa, John Pappajohn Business Building S362, Iowa
City, IA 52242-1994, United States. E-mail: \texttt{gao-zhengyuan@uiowa.edu}%
}%
\thanks{Southwestern University of Finance and Economics, RIEM Building 217,
Chengdu, 61000, P.R.China. E-mail: \texttt{zgao@swufe.edu.cn}%
}}
\maketitle
\begin{abstract}
This paper proposes a local representation for Empirical Likelihood (EL). 
EL admits the classical local linear
quadratic representation by its likelihood ratio property. A local estimator is
derived by using the new representation. Consistency, local asymptotic normality,
and asymptotic optimality results hold for the new estimator. In particular, when the regularity 
conditions do not include any differentiability assumption, these asymptotic results
are still valid for the local estimator. Simulations illustrate that the local method improves
the inference accuracy of EL.
\end{abstract}
\thispagestyle{empty}



\newpage{}

\baselineskip19.15pt \pagenumbering{arabic}

\section{Introduction }

A family of probability measures $\mathcal{E}_{\theta}=\{P_{\theta};\theta\in\Theta\}$
could represent a class of economic models. For a specific parameter
$\theta$ in $\Theta\in\mathbb{R}^{d}$, the probability $P_{\theta}$
measures the performance of the corresponding model. A sequence of
papers consider how to attain a suitable $P_{\theta}$ by comparing
a specified moment restriction function or moment constraint function
\[
\int m(x,\theta)dP_{\theta}(x)=\mathbb{E}_{\theta}[m(X,\theta)],
\]
to its sample counterpart\emph{
\[
\int m(x,\theta)dP_{n}(x)=\frac{1}{n}\sum_{i=1}^{n}m(X_{i},\theta),
\]
}where $P_{n}$ is the empirical distribution (empirical measure)
and $m(x,\theta)$ is a $k\times1$ vector with $k\geq d$ for given $x$ and $\theta$.%
\footnote{Although $P_{\theta}$ is indexed by $\theta$, the true distribution
of $m(X,\theta)$ does not depend on $\theta$. The notation $P_{\theta}$
can be interpreted as a pseudo measure of $m(X,\theta)$ and the specification
of this measure depends on the value of $\theta$. Later, $P_{\theta}$
is called an implied measure. %
} 

Empirical Likelihood (EL) fills in the gap between Generalized Methods
of Moments (GMM) and the classical Maximum Likelihood Estimation (MLE)
because it can incorporate the moment constraints into the classical
likelihood-based framework. \citet{QinLawless1994}, \citet{KitamuraStutzer1997},
and \citet{Smith1997} have shown that the estimators in both EL and
GMM-based estimates share many similar statistical features. As a
matter of fact, EL estimation with moment constraints has often been
recognized as a moment-based estimation method in econometrics. The
particular correspondence between $\mathcal{E}_{\theta}$ and $m(X,\theta)$
by EL is given as follows. For $n$ observations, the moment-based
EL is:
\[
\max_{\theta,p_{1},\dots,p_{n}}\left\{ \left.\prod_{i=1}^{n}np_{i}\right|\sum_{i=1}^{n}p_{i}m(X_{i},\theta)=0,\: p_{i}\geq0,\:\sum_{i=1}^{n}p_{i}=1\right\} .
\]
Function $m(X_{i},\theta)$ is of main interest in all moment-based
estimation methods. 

The connection between the moment-based estimation method and maximization
of likelihood ratios comes from dual parameters, that is, the parameters
in a dual problem. The dual problem in \citet{KitamuraStutzer1997}
shows an alternative way of incorporating moment constraints from
GMM. The moment constraints no longer appear directly in the objective
functions as in GMM or other minimum distance methods. The moment
constraints, however, are controlled dually by the Lagrangian multiplier
in EL and then appear indirectly in the modified objective functions.\footnote{Duality theory studies a\emph{ }pair of optimization problems, the
initial problem, which refers to the ``primal problem'',
and the dual problem. The aim of dual problem is to obtain more information
about the primal problem. For EL and its related methods, the information
of constraints and the information of optimal ``weights'' $\{p_{i}\}_{i\leq n}$
of these constraints are presented in a single criterion by the duality
theory.} Using the auxiliary dual parameters, \citet{Smith1997} and \citet{NeweySmith2004}
show that a class of estimators including Exponential Tilting, continuous
updating GMM and EL, will have better statistical properties than
the original GMM whose weighting matrices are not necessarily optimal.
However, the minimax type nonlinear optimization induced by the dual
parameters makes EL and its related methods less applicable. 

The main contributions of the paper are twofold. First, we present
a feasible local criterion EL function which resolves the minimax
criterion over nonlinear likelihood function. When the likelihood function
in the primal problem of EL has nonlinear constraints, the
objective function of the dual problem forms a minimax criterion with
an infinite dimensional (functional) dual parameter. Without an explicit
functional form, the dual parameter cannot be specifically incorporated
in a global representation. Furthermore, the dual parameter of EL may have unstable solution(s)
that give a thread to estimation and also a thread to
computation. Because the dual parameter appears in
the criterion function in the primal problem and also appears the Hessian matrix in the optimization
algorithm. The localization method will mitigate these threads. The basic idea in this paper
is to linearize the nonlinear optimization problem of EL by localizing
the likelihood ratio function. Once the nonlinear problem becomes
a linearized optimization problem, the minimax problem is reduced
to a linear or a quasi-linear programming problem.\footnote{In optimization, when one attempts to solve a nonlinear optimization problem, one should first think about transferring the problem into
a linear or quasi-linear environment.} 

The second contribution is to derive the local estimator, propose its computation
method and study its asymptotic properties.
The estimator comes from the primal-dual scheme together with Netwon-Le
Cam's localization. The estimation principle is as follows: approximate
the likelihood ratio in the primal problem, obtain a tractable dual
representation for the approximating primal problem, update the dual
parameter and then return its value to the primal problem. The dual
result follows the idea of the Kitamura-Stutzer \citep{KitamuraStutzer1997}
type duality and it assists to adjust the multiplier and the primal
likelihood function. Statistical properties of this iterative scheme
will depend only on the last iteration of the constructed estimator.
This estimator is asymptotically optimal. In addition, the local
estimator does not require a differentiable condition of the likelihood
function. This result could be important to practitioners. It provides
a theoretical ground for the practical use of EL estimator for data
with contaminated moment constraints which will be illustrated in Monte
Carlo simulations. 

In particular, localization representation avoids poor behaviors of
likelihood ratios in some corrupted models by contamination. In our
consideration, contamination induces non-informative likelihood ratio
values for estimation or poorly behaved Hessian matrices for computation.
For example, if the likelihood is flat in a neighborhood of some critical
points, the Hessian matrix is (near-) singular and the computation
may break down at these points.  In
the implementation, the likelihood of EL includes a vector of implied
probabilities $(\tilde{p}(X_{1},\theta),\dots,\tilde{p}(X_{n},\theta))$
where $\theta\in\Theta\subset\mathbb{R}^{d}$. Localization considers
the probability vector $(\tilde{p}(X_{1},\theta^{*}+\delta_{n}\tau),\dots,\tilde{p}(X_{n},\theta^{*}+\delta_{n}\tau))$
on a neighborhood of some $\theta^{*}$ and returns numbers for each
$\tau$ instead of functions. A well-behaved local representation
ensures the existence of the derivative of this representation. By
definition, when the derivative exists, small changes will not blow
up the approximation of the original likelihood ratio function and
this representation is therefore robust to these changes. Thus localization
avoids the peculiar points that break down the computational routines.

One could think of this local representation as an alternative
criterion function to the likelihood ratio. The following discusses
the connection between frequently used criterion functions and the
local approximating likelihood ratio criterion in this paper. EL has
been embedded into several general criteria, see e.g. \citet{Smith1997},
\citet{Baggerly1998}, \citet{NeweySmith2004}. The aims of these
estimation methods are similar: to optimize a criterion function of
$\theta$, such as a likelihood ratio function, subject to some constraint
of $m(X,\theta)$. The choice of criterion functions matters for the
efficiency and the robustness of an estimator. To balance the tradeoff
between these two objectives, \citet{Schennach2007} suggests a two-step
inference method by switching the empirical discrepancy between two
criterion functions, Kullback-Leibler and likelihood ratio. Although
this two-step inferential method works better than either its criterion
functions, changing the criterion function in the intermediate stage
could distort the supports of likelihood ratio and of Kubllback-Leibler
functions.\footnote{Kullback-Leibler and likelihood ratio use different measures as their
dominating measures in the criterion functions. Switching the position
of these measures require a mutual contiguity between the empirical
measure $P_{n}$ and implied probability measure $\tilde{P}_{n}$.
In other words, for every sequence $\{A_{n}\}_{n\in\mathbb{Z}}$,
$P_{n}(A_{n})\rightarrow0$ implies $\tilde{P}_{\theta}(A_{n})\rightarrow0$,
vice versa. This is a rather strong requirement even for a linear
constraint problem.} Instead of using two-step method, \citet{KitamuraOtsuEvdokimov2009}
suggest using Hellinger's distance as the criterion. Hellinger's distance
has a better topological structure than likelihood ratio and its estimator
shares almost the same first order statistical properties with EL.
In this paper, our representation of the classical likelihood ratio
is a linear-quadratic type approximation. This representation locally
obtains some Gaussian properties and therefore maintains a similar
topological structure as Hellinger's distance.\footnote{The covariance function of the approximating log-likelihood ratio
process can be attached to an inner product space (pre-Hilbert space) which is close to the $L^{2}$ structure considered by the Hellinger distance.}  The linear-quadratic representation induces
the Newton type iteration which is easier for implementations than
previous methods since it does not calculate the Hessian based on
the second derivative of moment constraints. 

The rest of the paper is organized as follows. Section 2 describes
EL and gives a version of consistency result without requiring the
existence of derivatives. Section 3 presents the local representation
of EL. Section 4 gives the local estimator and its asymptotic properties.
In Section 5 we describe two Monte Carlo experiments based on linear and nonlinear moment constraints.
Finally, conclusions appear in Section 6. Proofs are given in the
Appendix.

\section{Empirical Likelihood}
EL considers a finite dimensional parameter $\theta$ and an increasing
number of  \[
\mathbf{p}(X,\theta):=(p(X_{1},\theta),\dots,p(X_{n},\theta)).
\]
In this paper, the random variable $X_{i}$ is assumed to be i.i.d..
EL simultaneously finds the optimal $\theta$ and the optimal $\mathbf{p}(X,\theta)$
that satisfy the required moment constraints \[
\sum_{i=1}^{n}m(X_{i},\theta)p(X_{i},\theta)=0.
\] Its criterion is: \[
\sup_{p_{i},\theta}\left\{ \sum_{i=1}^{n}\log np_{i}\;|\; p_{i}\geq0,\quad\sum_{i=1}^{n}p_{i}=1,\sum_{i=1}^{n}p_{i}m_{i}(\theta)=0\right\} ,
\] where $p_{i}$ is a shorthand for $p(X_{i},\theta)$ given the value
$\theta$. An explicit expression for the optimal $p_{i}$'s can be
derived using the Lagrangian method and gives the solution:
\[
\tilde{p}_{i}(\theta):=\frac{1}{n}\frac{1}{1+\lambda_{n}^{T}m_{i}(\theta)},
\]
where $\tilde{p}_{i}(\theta)$ is called the \emph{implied probability}.
The candidate solutions belong to the family 
\[
\mathcal{E}_{\theta}:=\{\tilde{P}_{\theta}:\theta\in\Theta,\int m(X,\theta)d\tilde{P}_{\theta}=0,\tilde{P}_{\theta}\ll P_{0},\tilde{P}_{\theta}\ll P_{n}\},
\]
where $d\tilde{P}_{\theta}(x_{i})=\tilde{p}_{i}(\theta)d\mu$ for
a counting measure $\mu$.\footnote{The family $\mathcal{E}_{\theta}$ obtains both continuous measures
and discrete measures. The definition will become clear once we introduce
the infinite divisibility concept.%
} $\{\tilde{P}_{\theta}\ll P_{0},\tilde{P}_{\theta}\ll P_{n}\}$ means
that $\tilde{P}_{\theta}$ is contiguous with respect to both $P_{0}$
and $P_{n}$. For every sequence $\{A_{n}\}_{n\in\mathbb{Z}}$, $P_{n}(A_{n})\rightarrow0$
implies that $\tilde{P}_{\theta}(A_{n})\rightarrow0$ and meanwhile
$P_{0}(A_{n})\rightarrow0$ implies that $\tilde{P}_{\theta}(A_{n})\rightarrow0$.
The $\lambda_{n}$ is the solution of:
\begin{equation}
\frac{1}{n}\sum_{i=1}^{n}\left[\frac{m_{i}(\theta)}{1+\lambda_{n}^{T}m_{i}(\theta)}\right]=0.\label{eq:2-1-0}
\end{equation}
Let the average log-likelihood ratio of the implied probability between
any two parameter values $\theta_{1}$ and $\theta_{2}$ be:
\[
\Lambda_{n}(\theta_{1},\theta_{2}):=\frac{1}{n}\sum_{i=1}^{n}\log\left[\frac{\tilde{p}_{i}(\theta_{1})}{\tilde{p}_{i}(\theta_{2})}\right]
\]
and define the average log-likelihood ratio of the implied probability
given $\theta$ and counting numbers $1/n$ as 
\[
\Lambda_{n}(\theta):=\frac{1}{n}\sum_{i=1}^{n}\log n\tilde{p}_{i}(\theta).
\]
The constraint $0\leq\tilde{p}_{i}\leq1$ requires that the inequality
$1+\lambda_{n}^{T}m_{i}(\theta)\geq1/n$ always holds. The population
$\lambda(\theta):=\lim_{n\rightarrow\infty}\lambda_{n}$ must lie
in a convex and closed set $\Gamma_{\theta}=\lim_{n\rightarrow\infty}\cup_{i=1}^{n}\Gamma_{\theta,i}$.
For fixed $n$, the set $\Gamma_{\theta,n}$ is defined as a collection
of subsets of 
\[
\{\lambda_{n}:1+\lambda_{n}^{T}m(X_{i},\theta)\geq1/n,\: i=1,\dots,n,\:\theta\in\Theta\}.
\]

In the rest of this section, we derive an other consistency result
for EL estimation. Our intention is to obtain the consistency result
without assuming the differentiability of the moment restriction $m(X,\theta)$.
The differentiability is often assumed because it is a natural way
to derive an expansion of the objective function at the true parameter.
This expansion will link the asymptotic behaviors of $T_{n}-\theta_{0}$
with those of the sample averages of $\partial m(X,\theta)/\partial\theta$
and hence it is useful for proving both strongly and weakly convergences.
But as a trade-off, one needs to impose additional identification
conditions and limit distribution conditions for $n^{-1}\sum_{i}^{n}\partial m(X_{i},\theta)/\partial\theta$
and $n^{-\frac{1}{2}}\sum_{i}^{n}\partial m(X_{i},\theta)/\partial\theta$
respectively. Because our representation will not rely on such an
expansion, we weaken the conditions for consistency. 

There are many existing results of EL's consistency. \citet{KitamuraTripathiAhn2004}
relax the assumptions in \citet{QinLawless1994} and \citet{KitamuraStutzer1997}
and obtain consistency of the estimator based on Wald's approach \citep{Wald1949}.
\citet{NeweySmith2004} assume the differentiability of Lagrangian
multiplier rather than that of $m(X,\theta)$. However, due to the
non-analytical form of $\lambda(\theta)$, this assumption is quite
strong. \citet{Schennach2007} gives another consistency proof for
a non-differentiable objective function and avoids applications of
a Taylor expansion. The differentiability of the moment restriction,
however, is still assumed there in order to obtain a valid approximation
for the Lagrangian $\lambda(\theta)$. The conditions in the following
Theorem \ref{thm:consistency} are similar to the standard $M$-estimator
conditions in \citet{Huber1981}, thus the differentiability assumption
is not required. In order to ensure the EL estimator consistent for
this case, we need to give a result of EL consistency under weaker
conditions. Here are the conditions:
\begin{condition}
\label{con:1-1} 

$(i)$ $M(\theta):=\mathbb{E}[m(X,\theta)]$ exists for all $\theta\in\Theta$
and has a unique zero at $\theta=\theta_{0}.$ 

$(ii)$ $\theta_{0}$ is a well-separated point in $M(\theta)$ such
that
\[
\inf_{\theta:d(\theta,\theta_{0})\geq\epsilon}|M(\theta)|>|M(\theta_{0})|=0,
\]
where $\epsilon$ is an arbitrary value larger than zero and $d(\cdot,\cdot)$
is any distance function on $\Theta\times\Theta$. 

$(iii)$ $m(X,\theta)$ is continuous in $\theta,$ 
\[
\lim_{\theta'\rightarrow\theta}\left\Vert m(X,\theta)-m(X,\theta')\right\Vert =0.
\]

$(iv)$ Let $\infty$ be the one-point compactification of $\Theta$,
then there exists a continuous function $b(\theta)$ bounded away
from zero, such that $(1)$ $\sup_{\theta\in\Theta}m(X,\theta)/b(\theta)$ is integrable,  $(2)$ $\lim\inf_{\theta\rightarrow\infty}\left\Vert M(\theta)\right\Vert /b(\theta)$
is larger than $1$, and $(3)$ $\lim\sup_{\theta\rightarrow\infty}\left\Vert m(X,\theta)-M(\theta)\right\Vert /b(\theta)<1$.

$(v)$ $\sum_{i=1}^{n}\left[m(x_{i},\theta_{0})m(x_{i},\theta_{0})^{T}\right]/n$
is full rank for all $n\geq1$.
\end{condition}
Condition \ref{con:1-1} $(i)$ ensures the model is identified for
a small neighborhood of $\theta_{0}$. $(ii)$ is a local separability
condition. $(iii)$ is used to obtain the continuity of the Lagrangian
multiplier. $(iv)$ is an envelope assumption; it is used to obtain
some dominated convergence results. The one-point (Alexandroff) compactification
allows us to let $\theta$ approach any boundary place of $\Theta$,
even if $\Theta$ is not compact and may extend indefinitely. The
usual proof of EL consistency \citep{QinLawless1994} requires the
existence of the continuous derivative of $m(X,\theta)$ and that
the derivative is of full rank. Condition \ref{con:1-1} is less restrictive
because it allows for irregular cases where the usual ``delta method''
does not work, e.g. when $m(X,\theta)$ is non-differentiable. Condition
\ref{con:1-1} $(i)$-$(iv)$ are the standard M-estimator conditions
in \citet{Huber1981} and are very weak in the context of parametric
models. 
\begin{thm}
\label{thm:consistency}If Condition \ref{con:1-1} holds, then every
sequence $T_{n}$ satisfying 
\[
T_{n}:=\arg\sup_{\theta\in\Theta}\sum_{i=1}^{n}\log n\tilde{p}_{i}(\theta)=\arg\sup_{\theta\in\Theta}n\Lambda_{n}(\theta)
\]
 will converge to $\theta_{0}$ almost surely. 
\end{thm}
Note that this theorem does not require any differentiation condition.
However, the differentiability is implicitly obtained in the later
section. In fact, the ``local'' concept is the analog of ``differential''.
If one fixes a particular $\theta_{0}$ in $\Theta$ and investigates
what happens to the likelihood ratio function with parameter sequences
of the form $\theta=\theta_{0}+\delta_{n}\tau$, with $\delta_{n}\rightarrow0$
as $n$ goes to infinity, then $\delta_{n}$ yields a sort of differentiation
rate just as the differentiation rate in basic calculus, and then
the whole localization problem can be analyzed as a kind of differentiability
analysis for the likelihood ratio function. The term $\tau$ is called
local parameter since it is an index for local features. This technique
often appears in the evaluation of local power of test statistics
and statistical experiments, see \citet{Vaart1998} and \citet{LeCamYang2000}.

\section{Gaussian Properties and Localization of EL}

A non-closed form dual parameter $\lambda_{n}$ induces a non-closed form probability vector
$\tilde{\mathbf{p}}(\theta)$. General techniques such as empirical processes
of studying irregular behavior of the functions are also not directly
applicable because the functional form of $\lambda_{n}$ has no closed-form
representation, since it is the solution of Equation (\ref{eq:2-1-0})
that depends on the sample size and parameter values. In this section,
we propose alternative conditions and specifications of EL to standardize the problem.

\subsection{Approximation for an Infinitely Divisible Family}

Instead of studying the implied probability vectors $\tilde{\mathbf{p}}(\theta)$, we
consider a family of probability measures
\[
\mathcal{E}_{\theta}:=\{\tilde{P}_{\theta}:\theta\in\Theta,\int m(X,\theta)d\tilde{P}_{\theta}=0\},
\]
where the discrete vector $\tilde{\mathbf{p}}(\theta)$ satisfies
\[
\sum_{i}^{n}m(X_{i},\theta)\tilde{p}_{i}(\theta)=0.
\]
If a random variable $\xi$, for every natural number $n$, can be
represented as the sum
\[
\xi=\xi_{1,n}+\xi_{2,n}+\cdots+\xi_{n,n}
\]
of $n$ i.i.d random variables $\xi_{1,n},\dots,\xi_{n,n}$, then
$\xi$ is called \emph{infinitely divisible} \citep[p. 78]{GnedenkoKolmogorov1968}.
A probability distribution is said to be infinitely divisible if and
only if it can be represented as the distribution of the sum of an
arbitrary number of i.i.d random variables. A family of such distributions
is often referred to as an \emph{infinitely divisible} \emph{family}.
In our case, for arbitrary sample size $n$ and fixed $\theta$, the
log-likelihood ratio process is 
\[
\Lambda((X_{1},\dots,X_{n}),\theta)=\log n\tilde{p}(X_{1},\theta)+\cdots+\log n\tilde{p}(X_{n},\theta).
\]
Every additional term $\log n\tilde{p}(X_{i},\theta)$ is an identical
distributed increment of this log-likelihood ratio process. One crucial
deficiency of the above argument for EL is that $\tilde{p}(X_{i},\theta)$
are not independent for all $i$s. Because $\lambda_{n}$ appears in
$\tilde{p}(X_{i},\theta)$ for $i=1,\dots n$. But since the dependence is introduced by $\lambda_{n}$ only and $\lambda_{n}$ appears as the same form for all $\tilde{p}(X_{i},\theta)$, once the value of $\lambda_{n}$ is conditioning, the rest part of $\tilde{p}(X_{i},\theta)$ will be independent with $\tilde{p}(X_{j},\theta)$ for any $i\neq j$. 

For a sufficient large $n$ and a fixed $\theta$, $\lambda_{n}$
in $\log n\tilde{p}(X_{1},\theta),$ is a stochastic element.\footnote{In a localization approach, when $\theta$ is given, $\lambda_{n}$
will converge to a normal random variable with mean zero, see e.g.
Theorem 1 in \citet{QinLawless1994}.
} In this case, one can think that the integral of the log-likelihood
ratio process, $\sum_{i}\log n\tilde{p}(X_{i},\theta)$, represents
an infinite divisible process $\xi$ in $n$ additive terms $\xi_{1,n}+\xi_{2,n}+\cdots+\xi_{n,n}$.\footnote{More details about such a construction are discussed in \citet[Chapter 5]{LeCamYang2000},
although in most cases, they use $\log(1+(p_{\theta}/p_{\vartheta})^{-1/2}-1)$
instead of $\log(p_{\theta}/p_{\vartheta})$ directly.%
} Thus $\mathcal{E}_{\theta}$ does not merely include the family
of distributions that satisfy the constraint $\int m(x,\theta)d\tilde{P}_{\theta}(x)$,
it also requires the sample average of the log-likelihood ratio process
of $\tilde{P}_{\theta}$ to be infinitely divisible. It seems
that EL inherits the moment constraint from moment-based methods
and inherits the infinitely divisibility from likelihood ratio
based methods.

An infinitely divisible family $\mathcal{E}$ admits a representation
$\mathcal{E}=\mathcal{E}_{1}\times\cdots\times\mathcal{E}_{n}=\otimes_{i=1,\dots,n}\mathcal{E}_{i}$
based on $n$ copies of the so called divisor $\mathcal{E}_{i}$,
where $n$ could be arbitrarily large and $\times$ denotes the direct
product. The family $\mathcal{E}$ is called\emph{ divisible} with
divisor $\mathcal{E}_{i}$. There are several well known infinitely
divisible families, e.g. Poisson and Gaussian families. 

It has been proved by \citet[Theorem 17.5]{GnedenkoKolmogorov1968}
that any infinitely divisible family can be approximated by a finite
number of Poisson type measures. This result basically means that
the infinitely divisible family constructed by $\{\log n\tilde{p}(X,\theta)\}$
can be approximated by a finite number of Poisson measures.\footnote{We give a short description about Poissonization in the appendix.
Infinite divisible family holds for arbitary number of $n$, so the
approximation in principle should be valid for the finite many $n$.
} Poisson family relates to the Gaussian family via the Hellinger's
affinity. We will use this property to deduce a representation of
the likelihood ratio process. 
\begin{thm}
\label{thm:expression} If $\tilde{P}_{\theta}$ is infinitely divisible
then when $n\rightarrow\infty$, the log-likelihood $\log d\tilde{P}_{\theta+\delta_{n}\tau_{n}}/d\tilde{P}_{\theta}$
can be approximated by a linear quadratic expression such that the
difference 
\begin{equation}
\sum_{i=1}^{n}\log\frac{d\tilde{P}_{\theta+\delta_{n}\tau_{n}}}{d\tilde{P}_{\theta}}(x_{i})-\left[\tau_{n}^{T}S_{\theta,n}-\frac{1}{2}\tau_{n}^{T}K_{\theta,n}\tau_{n}\right]\label{eq:2-1-2}
\end{equation}
tends to zero in probability for any bounded sequence $\{\tau_{n}\}$
with a random vector $S_{\theta,n}$ and a deterministic matrix $K_{\theta,n}$. 
\end{thm}
The infinite divisible feature gives us a useful representation for
the likelihood ratio process, a linear quadratic expression with a
local parameter $\tau_{n}$. This representation is similar as the
linerization method based on Taylor's expansion, however it does not
require the differentiability of the implied probability. With this expression,
we can construct our estimator without bothering with non-linear optimization,
since the parameter in (\ref{eq:2-1-2}) is re-parametrized by $\tau_{n}$
which appears linearly and quadratically in the equation. Furthermore,
neither the computational algorithm nor the weakly convergent statistics
involve any differentiation requirements. 
\begin{rem}
In the proof, we will show a relation for univariate Gaussian families.
For any pair of Gaussian measures $G_{\theta}$ and $G_{\vartheta}$,
there will be a linear-quadratic expression to relate them. Therefore,
the integral of $(dG_{\theta}/dG_{\vartheta})^{1/2}$ w.r.t. $G_{\vartheta}$
will have a linear quadratic representation. Then we show that if
$\tilde{P}_{\theta}$ is infinitely divisible, $(d\tilde{P}_{\theta}/d\tilde{P}_{\vartheta})^{1/2}$
will be approximately equal to $(dG_{\theta}/dG_{\vartheta})^{1/2}$,
so $(d\tilde{P}_{\theta}/d\tilde{P}_{\vartheta})^{1/2}$ will also
have a linear quadratic representation.
\end{rem}

\begin{rem}
The linear-quadratic approximations to the log-likelihood ratios can
possibly be used with other minimum contrast estimators, but such
constructions only lead to asymptotically sufficient estimates, in
the sense of Le Cam, when the contrast function mimics the properties
of log-likelihood function, at least locally. 
\end{rem}

\begin{rem}
From a computational aspect, when confronted with the nonlinear optimization,
the Hessian matrix of the problem in some cases is difficult to evaluate
especially in regions that are either extremely flat or very erratic.
It is then computationally more efficient to consider the local optimization
and avoid a singular or non-invertible Hessian matrix rather than
calculate the global second order derivative of the objective function. 
\end{rem}

\begin{rem}
\label{Remark_covariance}Theorem \ref{thm:expression} shows that
with a proper choice of $\delta_{n}$, the log-likelihood ratio can
be approximated by a linear-quadratic representation. One of the main
focus of this representation is the quadratic term. For
a pair of Gaussian measures $(G_{\theta},G_{\vartheta})$ with dominating
measure $\mu$ we will have 
\allowdisplaybreaks
\begin{align}
 & \int\left(\frac{dG_{\theta}}{dG_{\vartheta}}\right)^{\frac{1}{2}}dG_{\vartheta}=\int dG_{\theta}^{\frac{1}{2}}dG_{\vartheta}^{\frac{1}{2}}d\mu\notag\\
= & \mathbb{E}\exp\left\{ \sum_{i=\theta,\vartheta}\frac{1}{2}\left[L(i)+\mathbb{E}\log(\frac{dG_{i}}{d\mu})\right]\right\} \notag\\
= & \left[\exp-\frac{1}{4}\left(K(\theta,\theta)+K(\vartheta,\vartheta)\right)\right]\cdot\mathbb{E}\exp\left(\sum_{i=\theta,\vartheta}\frac{1}{2}L(i)\right)\label{eq:step2-1}\\
= & \exp\left\{ \frac{1}{4}\left[2K(\theta,\vartheta)-K(\theta,\theta)-K(\vartheta,\vartheta)\right]\right\} ,\label{eq:step2-2}
\end{align}
where $L(i):=\{\log(dG_{i}/d\mu)-\mathbb{E}\log(dG_{i}/d\mu)\}$ for
$i=\theta,\vartheta$. The derivation of \eqref{eq:step2-2} is given in the Appendix. The property of $L(i)$ includes that it is Gaussian with expectation $\mathbb{E}L(i)=0$ and covariance kernel
$K(\theta,\vartheta)=\mathbb{E}L(\theta)L(\vartheta)$ and we have
$\mathbb{E}L(i)^{2}=K(i,i)$. Let
\[
q(\theta,\vartheta)=-8\log\int dG_{\theta}^{\frac{1}{2}}dG_{\vartheta}^{\frac{1}{2}}d\mu.
\]
Since the quadratic term is deterministic in the neighborhood of $\theta_{0}$,
we can use interpolation to find $K(\cdot,\cdot)$. With an arbitrary
mid-point $u$, three-point interpolation gives us:
\[
K(\theta,\vartheta)=-\left(q(\theta,\vartheta)-q(\theta,u)-q(u,\vartheta)\right).
\]
For small $|\theta-\vartheta|$, to speed up the computation, one
could use an approximated value $\Lambda_{n}(\theta,\vartheta)$ instead
of $q(\theta,\vartheta)$.\footnote{The concern is that the square root density computing may induce rounding
error. In fact $\frac{1}{2}\log\int(dG_{\theta}/dG_{\vartheta})^{1/2}dG_{\vartheta}$
approximately equals to $\frac{1}{2}\sum_{i}\log(dG_{\theta}/dG_{\vartheta})(x_{i})$
when $x_{i}$ is generated by $G_{\vartheta}$.}
\end{rem}

\subsection{Comparison with Other Conditions}

The standard EL ratio can be put into the form of the linear quadratic
representation in (\ref{eq:2-1-2}) but this requires some additional
assumptions, e.g. differentiability of $m(X,\theta)$. The following
proposition establishes this relation. 
\allowdisplaybreaks
\begin{prop}
\label{pro:expansion}Suppose that in addition to Condition \ref{con:1-1},
the following holds

(i) the model is just-identified, $\partial m(X,\theta)/\partial\theta<\infty$
for any $X$, the rank of $\mathbb{E}[\partial m(X,\theta)/\partial\theta]|_{\theta_{0}}$
equals $\dim(\theta),$ 

(ii) $\frac{1}{n}\sum_{i=1}^{n}[m_{i}(\theta)m_{i}(\theta)^{T}]$
and $\frac{1}{n}\sum_{i=1}^{n}[\lambda_{n}^{T}m_{i}(\theta)]^{2}$
are both finite for any positive $n$, even as $n\rightarrow\infty$,

then the log-likelihood ratio between $\tilde{p}_{\theta_{0}}$ and
$\tilde{p}_{\theta_{0}+\delta_{n}\tau}$ can be approximated by: 
\begin{align}
2\sum_{i=1}^{n}\log\frac{\tilde{p}_{\theta_{0}+\delta_{n}\tau_{n}}}{\tilde{p}_{\theta_{0}}}(x_{i})=\delta_{n}\tau_{n}^{T}A_{1} & +\frac{1}{2}\delta_{n}^{2}\tau_{n}^{T}A_{2}\tau_{n}^{T}+o_{p}(1)\label{eq:2-1}
\end{align}
where
\begin{align*}
A_{1} & =\mathbb{E}\frac{\partial m(X,\theta_{0})}{\partial\theta}^{T}\left(\mathbb{E}m(X,\theta_{0})m(X,\theta_{0})^{T}\right)^{-1}\sum_{i=1}^{n}m_{i}(\theta_{0}),\\
A_{2} & =\mathbb{E}\frac{\partial m(X,\theta_{0})}{\partial\theta}^{T}\left(\mathbb{E}m(X,\theta_{0})m(X,\theta_{0})^{T}\right)^{-1}\mathbb{E}\frac{\partial m(X,\theta_{0})}{\partial\theta^{T}}.
\end{align*}
\end{prop}
The expansion (\ref{eq:2-1}) is obtained simply by Taylor expansion
and the result therefore does not apply to the nonstandard problem
where the differentiability of $m(X,\theta)$ is questionable. However,
the result is intuitive as it mimics the standard Local Asymptotic
Normal (LAN) property for parametric models, see e.g. \citet[pp 104]{Vaart1998}.
The relation between (\ref{eq:2-1}) and (\ref{eq:2-1-2}) is also
quite clear: the first term is $\tau_{n}$ times a random vector,
and the second term is its variance. 
\begin{rem}
With the additional normality assumption on the average of $m_{i}(\theta_{0})$
and assuming $\delta_{n}=n^{-1/2}$ we will of course have:
\begin{align*}
\frac{1}{\sqrt{n}}\sum_{i=1}^{n} & \mathbb{E}\frac{\partial m(X,\theta_{0})}{\partial\theta}^{T}\left(\mathbb{E}m(X,\theta_{0})m(X,\theta_{0})^{T}\right)^{-1}m_{i}(\theta_{0})\\
\rightsquigarrow\mathcal{N} & \left(0,\mathbb{E}\frac{\partial m(X,\theta_{0})}{\partial\theta}^{T}\left(\mathbb{E}m(X,\theta_{0})m(X,\theta_{0})^{T}\right)^{-1}\mathbb{E}\frac{\partial m(X,\theta_{0})}{\partial\theta^{T}}\right).
\end{align*}
Asymptotic normality of the EL estimator is established by equation
(\ref{eq:2-1}) with additional conditions on the continuity or the
boundedness of second derivative of the moment restriction functions,
e.g. \citet{QinLawless1994}, \citet{NeweySmith2004} or \citet{KitamuraTripathiAhn2004}. 
\end{rem}

\begin{rem}
An alternative way of deducing this asymptotic normality is via Differentiability
in Quadratic Mean (DQM). This entails the existence of a vector of
measurable functions $S_{\theta_{0},n}$ such that
\begin{equation}
\int\left[\tilde{p}_{\theta_{0}+\delta_{n}\tau}^{1/2}-\tilde{p}_{\theta_{0}}^{1/2}-\frac{1}{2}\delta_{n}\tau^{T}S_{\theta_{0},n}\tilde{p}_{\theta_{0}}^{1/2}\right]^{2}d\mu=o(\Vert\delta_{n}\Vert^{2}),\label{eq:2-DQM}
\end{equation}
where $\delta_{n}\rightarrow0$. Note that the relation between the derivatives
of the square root density and the score function (when it exists) is:
\[
2\frac{1}{\sqrt{\tilde{p}_{\theta}}}\frac{\partial}{\partial\theta}\sqrt{\tilde{p}_{\theta}}=\frac{\partial}{\partial\theta}\log\tilde{p}_{\theta}.
\]
If along a path, the square root of the implied probability $\theta\mapsto\sqrt{\tilde{p}_{\theta}}$
is differentiable, then DQM basically means that a expansion of the
square root of $\tilde{p}_{\theta}$ is valid and the remainder term
is negligible in $L^{2}(\mu)$ norm. The term $S_{\theta,n}$ can
be considered as the score function of the implied probability $\tilde{p}_{\theta}$
at $\theta_{0}$. DQM implies that the condition does not require
the point-wise definition of the derivative of $m(\theta,X)$ therefore
it is less restrictive. 
\end{rem}
Suppose the implied probability includes the term $m(\theta,X)$ which
is not always differentiable. Then it deserves more efforts to relax
the restrictive condition on differentiability. In fact, Theorem \ref{thm:expression}
implies that the log-likelihood ratio belongs to the LAN family. The
result is already good enough for constructing an efficient (or asymptotic
sufficient) estimator. The expression in (\ref{eq:2-1-2}) is much
weaker than the regular conditions and DQM. It only states that log-likelihood
ratios of implied probabilities can be approximated by a linear-quadratic
expression.

\section{Local Estimation\label{sec:Local-Estimation}}

By the result (\ref{eq:2-1-2}) in Theorem \ref{thm:expression},
we can study the behavior of a pair $(\tilde{P}_{\theta+\delta_{n}\tau_{n}},\tilde{P}_{\theta})$
by looking at the log-likelihood ratio process $\Lambda_{n}(\theta+\delta_{n}\tau_{n},\theta)(X)$
with index $\tau_{n}$. The log-likelihood ratio process admits linear
quadratic approximations as $n\rightarrow\infty$, with the term $\tau_{n}S_{n}$
linear in $\tau_{n}$ and the term $\tau_{n}^{T}K_{n}\tau_{n}$ quadratic
in $\tau_{n}$. The numerical values of the approximation depend on
the concentrated point $\theta$ and its local neighborhoods. With
these ideas in mind, we will show the following steps of constructing
a local type estimator. The explanation of each step is given after
the definition. 
\begin{defn*}
Given Condition \ref{con:1-1}, we define the following Le Cam type
local EL estimator in $5$ steps:

Step 1. Find an auxiliary estimate $\theta_{n}^{*}$ using a \emph{$\delta_{n}$-consistent}
estimator and restricted such that it lies in $\Theta_{n}$ (a $\delta_{n}$-sparse
discretization of $\Theta$).

Step 2. Construct a matrix $K_{n}$ with $K_{n,i,j}=u_{i}^{T}K_{n}u_{j}$,
$i,j=1,2,\dots,d$, given by
\begin{align*}
K_{n,i,j}= & -\left\{ \Lambda_{n}[\theta_{n}^{*}+\delta_{n}(u_{i}+u_{j}),\theta_{n}^{*}]\right.\\
 & \left.-\Lambda_{n}[\theta_{n}^{*}+\delta_{n}u_{i},\theta_{n}^{*}]-\Lambda_{n}[\theta_{n}^{*}+\delta_{n}u_{j},\theta_{n}^{*}]\right\} 
\end{align*}
 and $\{u_{1},\dots,u_{d}\}$ is a set of directional vectors in $\mathbb{R}^{d}$.
$u_{i}$ is a step-size in selected in advance. 

Step 3. Construct the linear term:
\[
u_{j}^{T}S_{n}=\Lambda_{n}[\theta_{n}^{*}+\delta_{n}u_{j},\theta_{n}^{*}]+\frac{1}{2}K_{n,j,j}.
\]
 Since all the right hand side values are known, $S_{n}$ can be computed
and is a proper statistic.

Step 4. Construct the adjusted estimator:
\[
T_{n}=\theta_{n}^{*}+\delta_{n}K_{n}^{-1}S_{n}.
\]

\end{defn*}

\paragraph*{Step 1}

The \emph{$\delta_{n}$-sparse} (discretization of the) parameter
space in Step 1 is suggested by Le Cam (see \citet[p 125]{LeCamYang2000}).
It requires a sequence of subsets $\Theta_{n}\subset\Theta$ satisfying
the following conditions $(i)$ that for any $\theta\in\Theta$ and
any constant $b\in\mathbb{R}^{+}$, the ball $B(\theta,b\delta_{n})$
contains a finite number of elements of $\Theta_{n}$, independent
of $n$, and $(ii)$ that there exist a $c\in\mathbb{R}^{+}$ such
that any $\theta\in\Theta$ is within a distance $c\delta_{n}$ of
a point of $\Theta_{n}$. If we think of $\Theta_{n}$ as nodes of
a grid with a mesh that gets finer as $n$ increases, then $(i)$
says that the grid does not get too fine too fast and $(ii)$ says
that the mesh refines fast enough to have nodes close to any point
in the original space $\Theta$. In other words, asymptotically $\theta_{n}^{*}$
should be close enough to $\theta_{0}$. Another interpretation of
$\delta_{n}$-sparsity is from a Bayesian perspective. That is for
arbitrary priors, the corresponding posteriors essentially concentrate
on the small vicinities shrinking at the rate $\delta_{n}$.

\paragraph*{Step 2}

As in the Remark \ref{Remark_covariance}, the covariance matrix in
Step 2 is an analog to the covariance kernel in Gaussian processes.
For a stationary Gaussian process, the covariance kernel is smooth
and differentiable in quadratic mean, the covariance kernel can be
written as
\begin{align*}
\mbox{Cov} & \left(\frac{1}{\delta_{n}}(G_{\theta+u\delta_{n}}-G_{\theta}),\frac{1}{\delta_{n}}(G_{\vartheta+u\delta_{n}}-G_{\vartheta})\right)\\
= & \frac{1}{\delta_{n}^{2}}\left(2C(\theta-\vartheta)-C(\theta-\vartheta+u\delta_{n})-C(\theta-\vartheta-u\delta_{n})\right)\\
\rightarrow & -\left.\frac{\partial^{2}C(h)}{\partial h^{2}}\right|_{h=\theta-\vartheta},
\end{align*}
where $C(\theta,\vartheta):=\mbox{Cov}\{G_{\theta},G_{\vartheta}\}$.
Since $K_{n}$ is an analog to the covariance kernel, the construction
of $K_{n}$ is nothing else but a finite difference of $\Lambda_{n}(\cdot,\cdot)$
which is analogous to the second derivative of the covariance kernel.

\paragraph*{Step 3 and 4}

With a control term $K_{n}$ which is asymptotically determined, all
the randomness of the log-likelihood ratio is contained in the first
term, $S_{n}$. Step 3 is to extract the randomness from $\Lambda_{n}(\cdot,\cdot)$
and construct the linear term. Step 4 is to construct the estimator.
To verify these two steps, we need to ensure that the covariance kernel
in (\ref{eq:2-1-2}) is invertible. 
\begin{prop}
\label{pro:LAQ(iii)}The matrices $K_{\theta,n}$ in (\ref{eq:2-1-2})
are almost surely positive definite. Any cluster point $K_{\theta}$
of $K_{\theta,n}$ in $P_{\theta,n}$-law is invertible.
\end{prop}
If $K_{n}-K_{\theta,n}$ converges to zero, then $K_{n}$ is also
invertible. This result will be given in the following Theorem \ref{thm:localEL}.
If $K_{n}$ is positive definite, by substituting $S_{n}=K_{n}\delta_{n}^{-1}(T_{n}-\theta_{n}^{*})$
into the linear quadratic expression: 
\begin{align*}
\tau_{n}^{T}K_{n}\delta_{n}^{-1}(T_{n}-\theta_{n}^{*})-\frac{1}{2}\tau_{n}^{T}K_{n}\tau_{n}= & -\frac{1}{2}\delta_{n}^{-2}\left[T_{n}-(\theta_{n}^{*}+\delta_{n}\tau_{n})\right]^{T}K_{n}\times\\
\left[T_{n}-(\theta_{n}^{*}+\delta_{n}\tau_{n})\right] & +\frac{1}{2}\delta_{n}^{-2}\left[T_{n}-\theta_{n}^{*}\right]^{T}K_{n}\left[T_{n}-\theta_{n}^{*}\right],
\end{align*}
we have a quadratic expression of $T_{n}$ and $(\theta_{n}^{*}+\delta_{n}\tau_{n})$.
The maximal value of this approximating representation of the log-likelihood
ratio is achieved when $\theta_{n}^{*}+\delta_{n}\tau_{n}=T_{n}$.
In other words, $\delta_{n}^{-1}(T_{n}-\theta_{n}^{*})$ is the estimator
for the local parameter $\tau_{n}$. 
\begin{rem}
The construction was originally proposed by \citet{LeCam1974}. He
supposed that there is a special interest in the likelihood function
at particular points where Taylor's expansion fails, e.g. for the
Laplace distribution. The advantage of the construction is that the
quadratic term does not depend very much on the particular auxiliary
estimation method that is used to obtain the value of $\theta_{n}^{*}$
and the construction is only determined in a local neighborhood of
the particular point.
\end{rem}

\begin{rem}
One may be concerned with the $\delta_{n}$-consistency requirement
for the auxiliary estimator. For a simple i.i.d. case, the $\delta_{n}$
is set to $n^{-1/2}$, the requirement is the same as asking for an
$\sqrt{n}$-consistent auxiliary estimator. Any $\sqrt{n}$-consistent
estimator should be, in principle, good enough from the estimation
perspective, because the auxiliary estimator $\theta_{n}^{*}$ is
at least in a neighborhood of $\theta_{0}$. However, in practice,
it may be hard to find a well behaved moment restriction function
around $\theta_{0}$. The use of local EL estimator is to overcomes
the problem and improve the auxiliary estimator. We suppose that $\theta_{n}^{*}$
is located within a range $n^{-1/2}$ of the true value, then a local
method would give a refinement. When consistency and asymptotic normality
are treated separately, one could take good care of consistency first
and then use localization method to improve the final result or one
could take care of the concentration of distribution first and then
correct the bias by localization.\end{rem}
\begin{thm}
\label{thm:localEL}Given Condition \ref{con:1-1}, $T_{n}$, $S_{n}$
and $K_{n}$ have following properties:

(i) $K_{n}^{-1}S_{n}-K_{\theta,n}^{-1}S_{\theta,n}$ and $K_{n}-K_{\theta,n}$
converge to zero in $\tilde{P}_{\theta,n}$-law where $(K_{\theta,n},S_{\theta,n})$
is in (\ref{eq:2-1-2}).

(ii) $\delta_{n}^{-1}(T_{n}-\theta)$ is bounded in $\tilde{P}_{\theta,n}$
-law.

(iii) if Equation (\ref{eq:2-DQM}) holds and the moment restrictions
are just-identifying, the sequence of models $\{\tilde{P}_{\theta,n}:\theta\in\Theta\}$
is LAN and 
\[
\delta_{n}^{-1}(T_{n}-\theta_{0})\rightsquigarrow\mathcal{N}(0,\Omega)
\]
 where $\Omega=\mathbb{E}\frac{\partial m(x,\theta_{0})}{\partial\theta}^{T}\left(\mathbb{E}m(x,\theta_{0})m(x,\theta_{0})^{T}\right)^{-1}\mathbb{E}\frac{\partial m(x,\theta_{0})}{\partial\theta^{T}}$. 
\end{thm}
The LAN theory is useful in showing that many statistical models can
be approximated by Gaussian models. In the parametric likelihood framework,
when the original model $P_{\theta}$ is smooth in the parameters,
i.e. DQM, the local parameter $\tau_{n}=\delta_{n}^{-1}(\theta_{0}-\theta_{n}^{*})$
can be used to construct a log likelihood ratio based on $P_{\theta_{0}+\tau_{n}\delta_{n}}$
that is asymptotically $\mathcal{N}(\tau_{n},I_{\theta_{0}}^{-1})$.
Here we use LAN in a moment based setting without further parametric
assumptions. Once LAN is established, asymptotic optimality of estimators
and of tests can be expressed in terms of LAN properties. 
\begin{rem}
Some other articles also utilize local information based on an EL
framework. \citet{DonaldImbensNewey2003} propose resampling data
from a local EL estimated distribution. \citet{KitamuraTripathiAhn2004}
consider another localized EL based on conditional moment restrictions
and use them to re-construct a smooth global profile likelihood function.
\citet{Smith2005} extends moment smoothing to GEL. These methods
construct smooth objective functions, implicitly or explicitly. Our
solution is to discretize the parameter space and then construct local
log-likelihood ratios as local objective functions. Thus our localization is viewing 
a different aspect of the problem. \end{rem}
Theorem \ref{thm:localEL} gives an asymptotic result on the weak
convergence of the estimator. In the theorem, the limit distribution
is based on a kind of Cram\'{e}r-Rao type lower bound and is essentially
a point-wise result. In order to obtain a result in a neighborhood
rather than at a single point, we will now state and prove a minimax
type theorem on the risk of any estimator.

Before giving the theorem, we need to introduce a technical concept
of \emph{$\delta_{n}$-regularity}. This concept expresses the desirable
requirement that a small change in the parameter should not change
the distribution of estimator too much. For the estimator sequence
$T_{n}$, if the difference between the distributions of $\delta_{n}^{-1}(T_{n}-\theta_{0}-\delta_{n}\tau)$
and $\delta_{n}^{-1}(T_{n}-\theta_{0})$ tends to zero under $P_{\theta_{0}+\delta_{n}\tau,n}$-law
and $P_{\theta_{0},n}$-law respectively, then $T_{n}$ is called
$\delta_{n}$-regular at the point $\theta_{0}$. 
\begin{thm}
\label{thm:minimax}Given Condition \ref{con:1-1} and letting $W$
be a non-negative bowl shaped loss function, if $T_{n}$ is $\delta_{n}$-regular
on all $\Theta$, then for any estimator sequence $Z_{n}$ of $\tau$,
one has 
\[
\lim_{b\rightarrow\infty}\lim_{c\rightarrow\infty}\lim_{n\rightarrow\infty}\inf_{n}\sup_{|\tau|\leq c}\mathbb{E}_{\theta_{0}+\delta_{n}\tau}\left[\min(b,W(Z_{n}-\tau))\right]\geq\mathbb{E}[W(\xi)]
\]
 where $\xi$ has a Gaussian distribution $\mathcal{N}(0,K^{-1})$.
The lower bound is achieved by $Z_{n}=\delta_{n}^{-1}(T_{n}-\theta_{0})$. 
\end{thm}
A loss function is ``bowl-shaped'' if the sublevel sets $\{u:W(u)\leq a\}$
are convex and symmetric around the origin. The value $b$ is used
to construct a bounded function $\min(b,W(Z_{n}-\tau_{n}))$. We let
$c$ go to infinity in order to cover a general case. The expectation
$\mathbb{E}_{\theta_{0}+\delta_{n}\tau}[\cdot]$ is taken w.r.t. a
measure $\mathcal{M}$ of the set $\{\theta:|\theta-\theta_{0}|\leq\delta_{n}\tau_{n}\}$
while $\mathbb{E}[\cdot]$ is taken w.r.t. a distribution of $K^{-1/2}\times\mathcal{N}(0,I)$
on $\xi$. 

The theorem can be interpreted as follows. When using the auxiliary
estimator $\theta_{n}^{*}$ in the likelihood ratio, this induces
randomness to the local parameter $\tau_{n}$. By using the LAN result
in Theorem \ref{thm:localEL}, we can attach the local parameter $\tau_{n}$
with a Gaussian measure. By the Gaussian prior assumption of $\tau_{n}$,
one can express the convergent procedure as a procedure of updating
a Gaussian prior, while for a centered Gaussian prior, this procedure
is to update the prior covariance matrix $\Gamma^{-1}$. The $\delta_{n}$-regularity
condition implies that $K_{n}$ will converge uniformly in a neighborhood
of $\theta_{0}$ for arbitrary measure $\mathcal{M}$. Thus the covariance
will converge to a the posterior covariance matrix $(K+\Gamma)^{-1}$.
The Gaussian randomness introduces a new random variable $\xi$ that
has the posterior covariance matrix $(K+\Gamma)^{-1}$. The lower
bound of the Bayes risk of this Gaussian variable is obtained by letting
$\Gamma$ go to zero, corresponding to initial values of $\tau$ widely
spread. This is the local asymptotic minimax theorem. It is based
on the minimax criterion and gives a lower bound for the maximum risk
over a small neighborhood of the parameter $\theta$. Because the
local EL can achieve this lower bound, it is an asymptotically optimal
estimator.

\section{Simulations}

\begin{table}
\vspace{6pt}

\caption{\label{tab:Iteration}}

\vspace{6pt}

\centering{}{\scriptsize{}}%
\begin{tabular}{cccc}
\multicolumn{4}{c}{{\footnotesize{Local Iteration}}}\tabularnewline
\hline 
\hline 
{\scriptsize{Iter. num.}} & {\scriptsize{$\lambda$}} & {\scriptsize{$\tau$}} & {\scriptsize{Local Estimator}}\tabularnewline
\hline 
{\scriptsize{1}} & {\scriptsize{0.148899}} & {\scriptsize{0.128070}} & {\scriptsize{2.069817}}\tabularnewline
{\scriptsize{2}} & {\scriptsize{0.134007}} & {\scriptsize{0.115263}} & {\scriptsize{2.058290}}\tabularnewline
{\scriptsize{3}} & {\scriptsize{0.120464}} & {\scriptsize{0.103737}} & {\scriptsize{2.047917}}\tabularnewline
$\vdots$ & $\vdots$ & $\vdots$ & $\vdots$\tabularnewline
{\scriptsize{12}} & {\scriptsize{0.053311}} & {\scriptsize{0.084166}} & {\scriptsize{1.992690}}\tabularnewline
{\scriptsize{13}} & {\scriptsize{0.043537}} & {\scriptsize{0.000000}} & {\scriptsize{1.992690}}\tabularnewline
\hline 
\end{tabular}
\end{table}

Throughout the paper, our concerns are the violations of the standard
regularity conditions for the moment restriction functions and their
derivative functions. In this section, we simulate two models whose
moment conditions are contaminated by some outliers. We call these
models contaminated models. The contamination in this experiment occurs
at a certain probability no matter how large the sample size is. The
simulations try to mimic the environment that few observations may
violate the boundedness condition for $m(X,\theta)$ and such observations
are not caused by the small number of samples. In other words, some
large values of $m_{i}(\theta)$ are caused by some $x_{i}$s and
these $x_{i}$s are systematically existing.

These features imply that the specification of the constraint $\mathbb{E}[m(X,\theta)]=0$
is invalid for the whole sample, although the specification is valid
for the uncontaminated sample. A completely misspecified model is
not of our interest. In our experiment, the contamination level is
controlled to a small value so that the model is not significantly
misspecified. A consequence of the mildly misspecified constraint
is that the moment-based estimators are biased. 

The full description of the localized EL's implementation is given
in the Appendix. From each iteration in the localization steps, the
value of the local estimator is adjusted. Table \ref{tab:Iteration}
gives an example of the information used in the localization step,
where $\lambda$ is the Lagrangian multiplier and $\tau$ is the local
parameter. The true value of the parameter is $2$. Due the mildly
mis-speciation of the moment restriction, $\lambda_{n}$ does not
reach $0$ when the estimator converges to the true value. However,
the local iteration of $\tau$ induce an almost unbiased estimate
result with the maximum local likelihood. 

\begin{figure}
\begin{centering}
\includegraphics[width=0.8\textwidth]{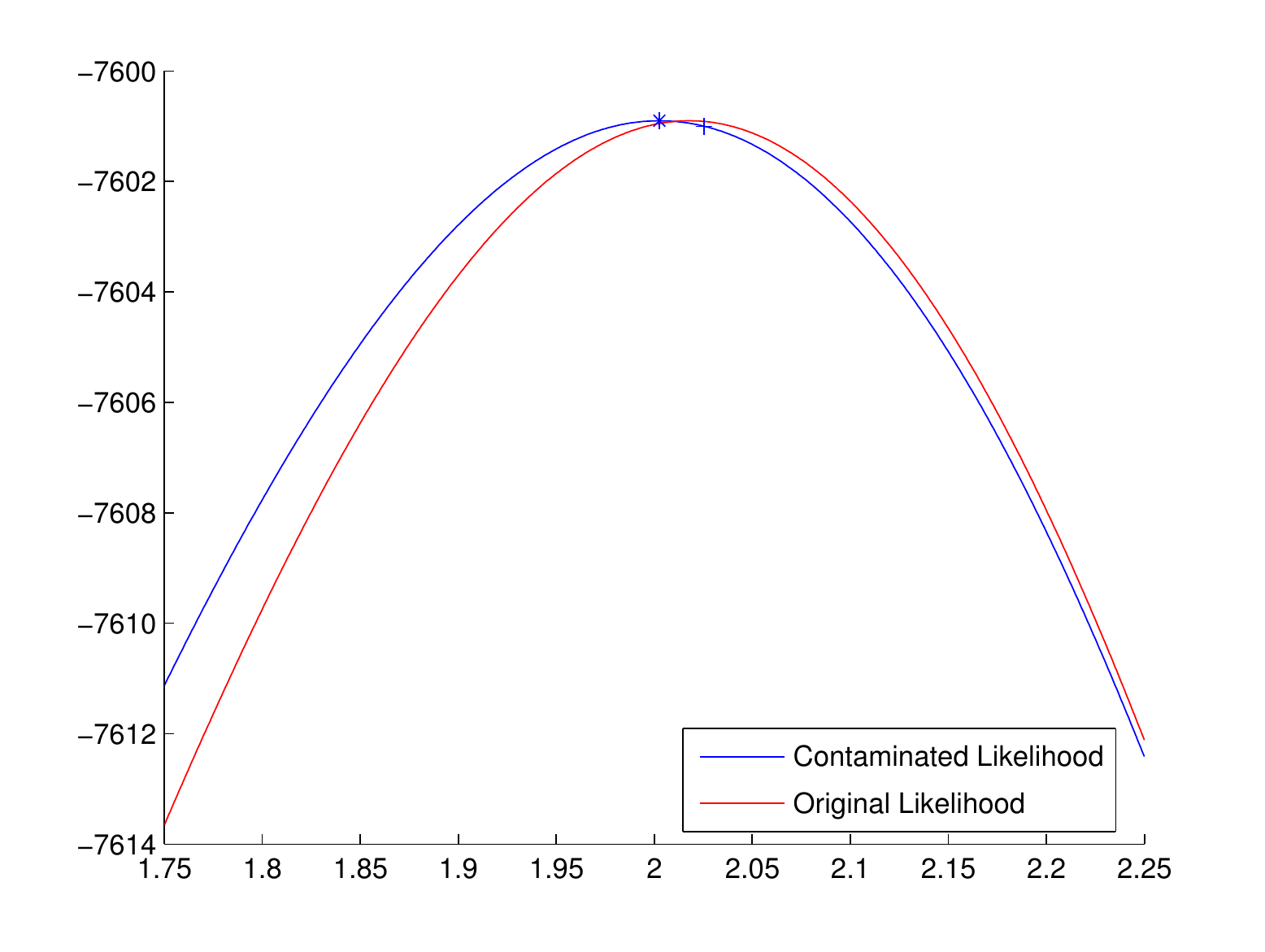}
\par\end{centering}

\begin{centering}
\includegraphics[width=0.8\textwidth]{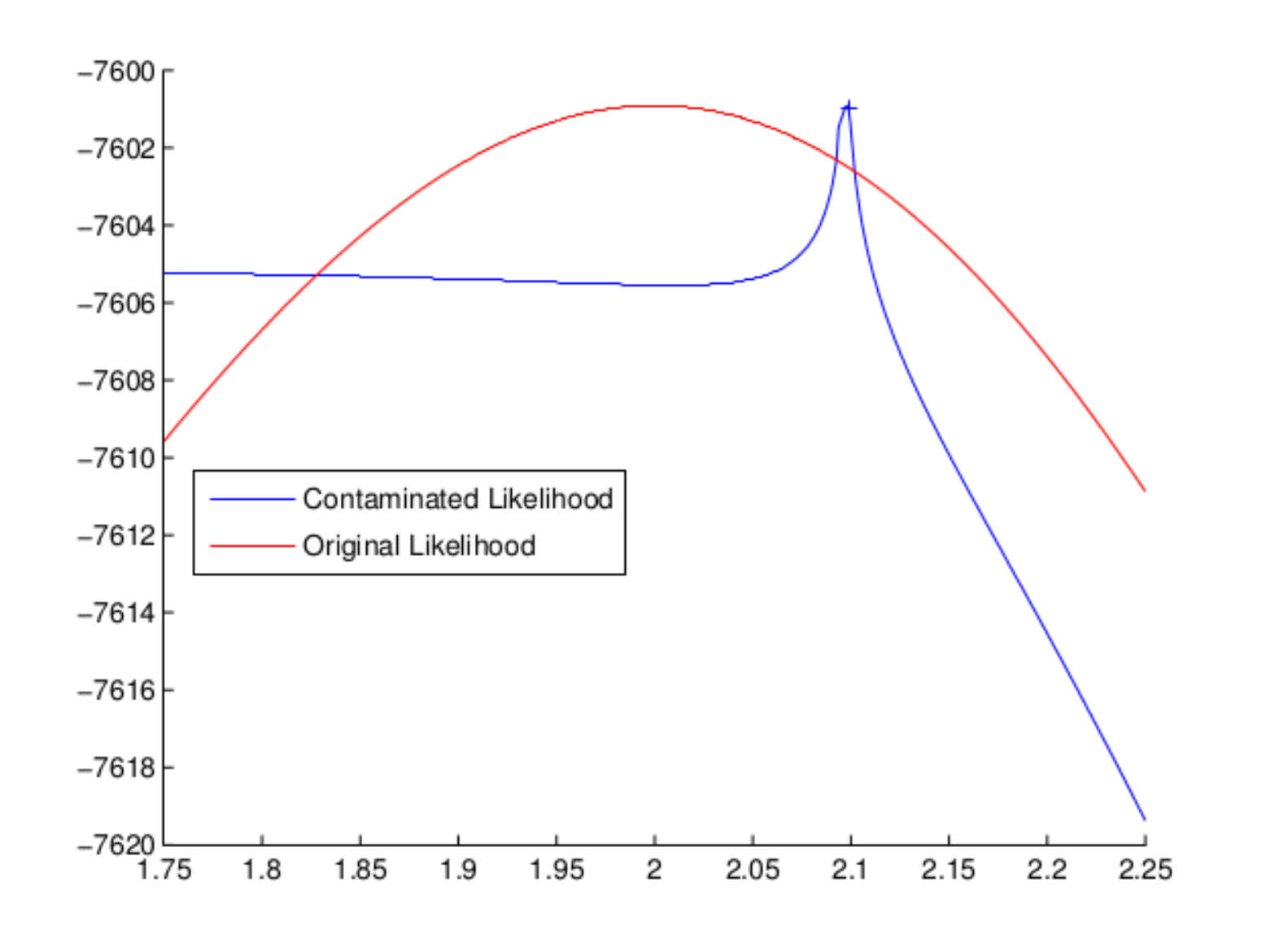}
\par\end{centering}

\caption{\label{fig:likelihoodplot}Log-likelihood. Cross stands for the LEL
estimation result and star stands for the EL estimation result. Blue
(Red) line is the log-likelihood for a contaminated (uncontaminated)
sample.}
\end{figure}

Figure \ref{fig:likelihoodplot} gives two representative phenomena
in the numerical experiment. When the simulation does not induce a
peculiar optimal point of the log-likelihood, EL rather than local
EL reach the peak of the empirical log-likelihood function. However,
such a peak is for the contaminated sample which induces misspecified
moment restrictions. This peak does not lead to the best solution.
Another situation is for irregular log-likelihood shape. In this case
the EL estimation does not give a local optimal answer, nor even report
a correct log-likelihood value. The problem is caused by the irregular
shape of the likelihood. The flat log-likelihood region and the non-smooth
peak break down the global search routine in the EL estimation. Although
local EL estimate value does not correspond to the parameter value
that gives the optimal log-likelihood for uncontaminated sample, local
EL estimator reaches the local optimal point of the empirical log-likelihood
function.

\subsection{A Linear Experiment}

\begin{table*}
\vspace{6pt}

\caption{\label{tab:Estimation-results}}

$\:$\vspace{6pt}

\begin{centering}
{\scriptsize{}}%
\begin{tabular}{ccccccccc}
\multicolumn{9}{c}{{\footnotesize{Linear Model: LS Auxiliary Estimator}}}\tabularnewline
\hline 
\hline 
 & \multicolumn{4}{c}{{\scriptsize{c\%=0.5\%, L=10 (case I)}}} & \multicolumn{4}{c}{{\scriptsize{c\%=0.005\%, L=10000 (case II)}}}\tabularnewline
\hline 
{\scriptsize{Method}} & {\scriptsize{Mean}} & {\scriptsize{Median}} & {\scriptsize{MSE}} & {\scriptsize{IQR}} & {\scriptsize{Mean}} & {\scriptsize{Median}} & {\scriptsize{MSE}} & {\scriptsize{IQR}}\tabularnewline
{\scriptsize{LS}} & {\scriptsize{2.092564}} & {\scriptsize{2.092765}} & {\scriptsize{0.009234}} & {\scriptsize{0.033984}} & {\scriptsize{2.091643}} & {\scriptsize{2.094945}} & {\scriptsize{0.009237}} & {\scriptsize{0.036843}}\tabularnewline
{\scriptsize{IV}} & {\scriptsize{2.012655}} & {\scriptsize{2.012436}} & {\scriptsize{0.008073}} & {\scriptsize{0.130901}} & {\scriptsize{2.008087}} & {\scriptsize{2.012734}} & {\scriptsize{0.009991}} & {\scriptsize{0.148564}}\tabularnewline
{\scriptsize{EL}} & {\scriptsize{1.990643}} & {\scriptsize{2.011329}} & {\scriptsize{0.024535}} & {\scriptsize{0.128907}} & {\scriptsize{1.984457}} & {\scriptsize{2.008713}} & {\scriptsize{0.033921}} & {\scriptsize{0.146062}}\tabularnewline
{\scriptsize{Local EL(LS)}} & {\scriptsize{2.045564}} & {\scriptsize{2.050873}} & {\scriptsize{0.005683}} & {\scriptsize{0.077868}} & {\scriptsize{2.048954}} & {\scriptsize{2.056338}} & {\scriptsize{0.006055}} & {\scriptsize{0.084951}}\tabularnewline
\hline 
\end{tabular}
\par\end{centering}{\scriptsize \par}

\vspace{6pt}

\begin{centering}
{\scriptsize{}}%
\begin{tabular}{ccccccccc}
\multicolumn{9}{c}{{\footnotesize{Linear Model: IV Auxiliary Estimator}}}\tabularnewline
\hline 
\hline 
 & \multicolumn{4}{c}{{\scriptsize{c\%=0.05\%, L=100 (case III)}}} & \multicolumn{4}{c}{{\scriptsize{c\%=0.01\%, L=10000 (case IV)}}}\tabularnewline
\hline 
{\scriptsize{Method}} & {\scriptsize{Mean}} & {\scriptsize{Median}} & {\scriptsize{MSE}} & {\scriptsize{IQR}} & {\scriptsize{Mean}} & {\scriptsize{Median}} & {\scriptsize{MSE}} & {\scriptsize{IQR}}\tabularnewline
{\scriptsize{LS}} & {\scriptsize{2.091583}} & {\scriptsize{2.092988}} & {\scriptsize{0.009233}} & {\scriptsize{0.038327}} & {\scriptsize{2.091021}} & {\scriptsize{2.091687}} & {\scriptsize{0.009196}} & {\scriptsize{0.040935}}\tabularnewline
{\scriptsize{IV}} & {\scriptsize{2.001433}} & {\scriptsize{2.003463}} & {\scriptsize{0.009768}} & {\scriptsize{0.132206}} & {\scriptsize{2.008599}} & {\scriptsize{2.001776}} & {\scriptsize{0.010176}} & {\scriptsize{0.137512}}\tabularnewline
{\scriptsize{EL}} & {\scriptsize{1.985619}} & {\scriptsize{2.002222}} & {\scriptsize{0.028999}} & {\scriptsize{0.135985}} & {\scriptsize{1.983376}} & {\scriptsize{2.000420}} & {\scriptsize{0.032342}} & {\scriptsize{0.141040}}\tabularnewline
{\scriptsize{Local EL(IV)}} & {\scriptsize{2.011766}} & {\scriptsize{2.015890}} & {\scriptsize{0.008267}} & {\scriptsize{0.120573}} & {\scriptsize{2.026002}} & {\scriptsize{2.021211}} & {\scriptsize{0.006859}} & {\scriptsize{0.105591}}\tabularnewline
\hline 
\end{tabular}
\par\end{centering}{\scriptsize \par}

$\ $
\end{table*}

We consider a simple structural model with a $n\times1$ explanatory
vector $\mathbf{x}_{n}=(x_{1},\dots,x_{n})^{T}$, a $n\times1$ instrument
vector $\mathbf{z}_{n}$ and a disturbance vector $\mathbf{u}_{n}$,
$n=1000$. The parameter $\theta$ is equal to $2$. The $n\times1$
random vector $\varepsilon$ is assumed to be normal. The model is
as follows:
\begin{align*}
y_{i} & =x_{i}\theta+\varepsilon_{i},\\
x_{i} & =z_{i}\pi+u_{i}.
\end{align*}
In our numerical experiment, $\pi$ is set to one. The instrument
$\mathbf{z}_{n}$ is a design vector with a constant vector plus a
small noise and $\mathbf{z}_{n}$ is independent of $\varepsilon$
and $\mathbf{u}_{n}$. The uncertainty vector $\mathbf{u}_{n}$ is
a mixture of two normally distributed vector $\mathbf{u}_{n}^{(1)}$
and $\mathbf{u}_{n}^{(2)}$ where $u_{i}^{(1)}\sim\mathcal{N}(0,1)$
and $u_{i}^{(2)}\sim\mathcal{N}(L,1)$. $L$ is referred to the degree
of contamination. We introduce $u_{i}^{(2)}$ to generate a mis-specified
moment. In this experiment, $\mathbf{u}_{n}^{(2)}$ is a contaminated
element. The mixing rate of $\mathbf{u}_{n}^{(2)}$ in $\mathbf{u}_{n}$
is the probability of contamination. Let $c$ denote this probability.
If $P_{\mathbf{u}_{n}^{(i)}}$ denotes the distribution of $\mathbf{u}_{n}^{(i)}$,
then $P_{\mathbf{u}_{n}}=(1-c)P_{\mathbf{u}_{n}^{(1)}}+cP_{\mathbf{u}_{n}^{(2)}}.$
We impose the correlation between $\varepsilon$ and $\mathbf{u}_{n}$
by using the equation $\varepsilon_{n}=R\times\mathbf{u}_{n}+\varepsilon_{n}^{'}$
where $\varepsilon_{n}^{'}\sim\mathcal{N}(0,1)$ is independent of
$\mathbf{u}_{n}$. The covariance value $R$ is set to $0.1$.

The moment restriction function in this example is $\mathbf{z}_{n}^{T}(\mathbf{y}_{n}-\mathbf{x}_{n}\theta)$.
We will consider four different estimation methods, Least Squares
(LS), Instrumental Variables (IV), EL, and local EL.\footnote{In this setup, the IV estimator asymptotically has a degenerated second
moment. Thus in order to make a fair comparison, we only consider
the cases where the IV estimators are not widely spreaded.%
}  The estimators for LS, IV, EL are respectively
$(\mathbf{x}_{n}^{T}\mathbf{x}_{n})^{-1}\mathbf{x}_{n}^{T}\mathbf{y}_{n}$,
$(\mathbf{z}_{n}^{T}\mathbf{x}_{n})^{-1}\mathbf{z}_{n}^{T}\mathbf{y}_{n}$
and \[ \min_{\beta}\max_{\lambda_{n}}\sum_{i=1}^{n}\log(1+\lambda_{n}m_{i}(\theta)),\] where $m_{i}(\theta)=z_{i}(y_{i}-x_{i}\theta)$. 

\begin{figure}
\begin{centering}
\includegraphics[width=0.8\textwidth]{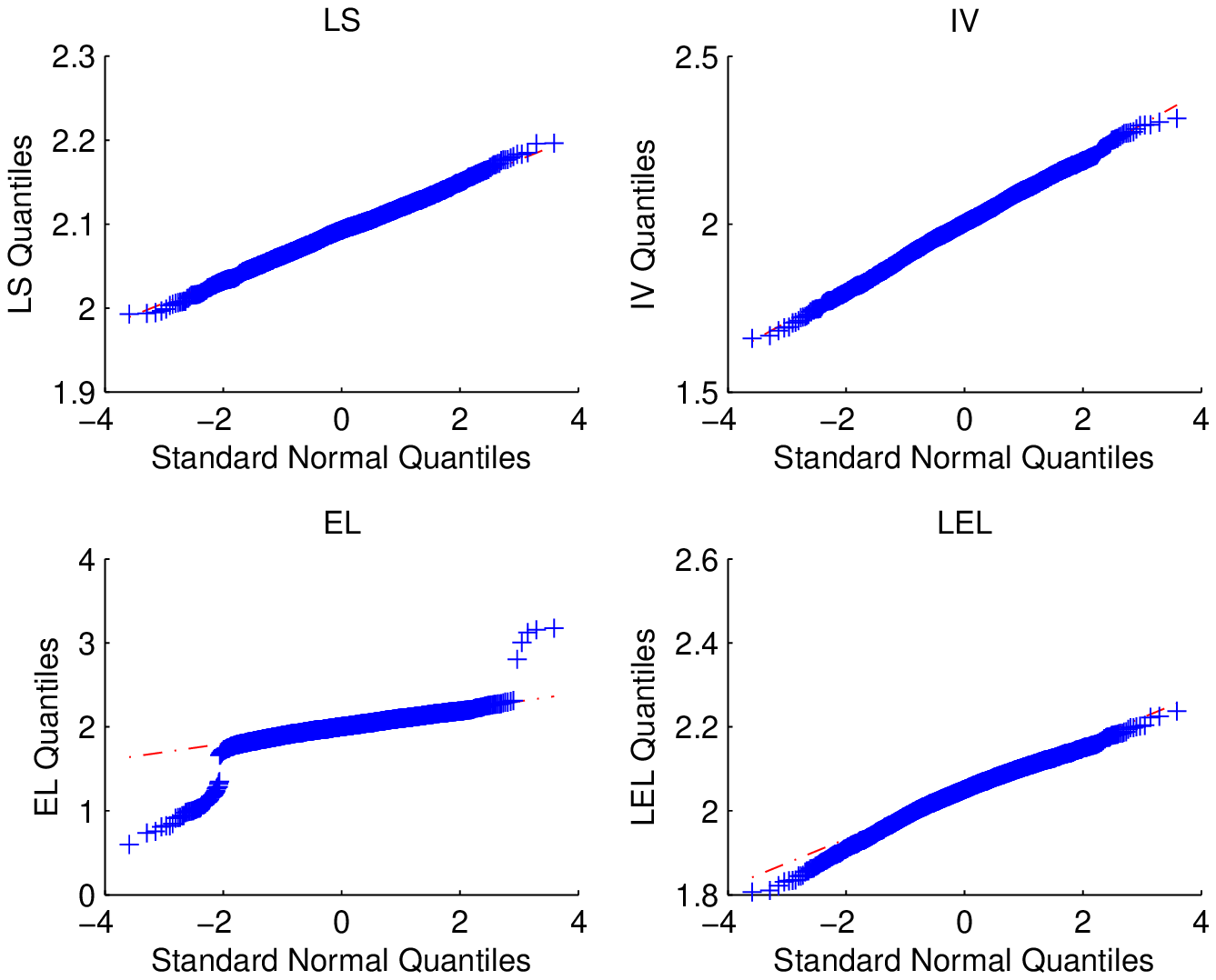}
\par\end{centering}

\begin{centering}
\includegraphics[width=0.8\textwidth]{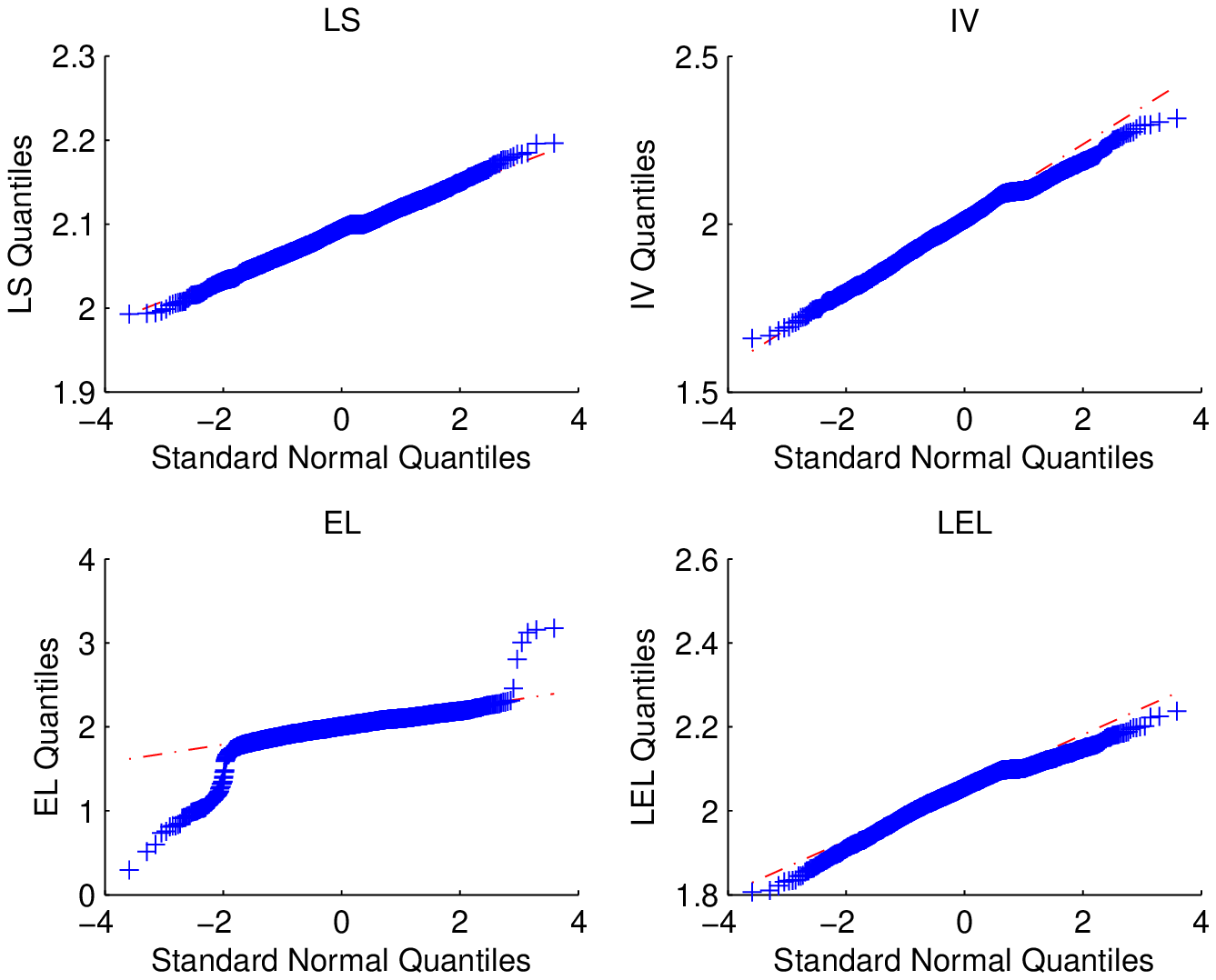}
\par\end{centering}

\caption{\label{fig:qqplot-ls}QQ plot (Densities of estimators). Up: case
(I). Down: case (II).}
\end{figure}

\begin{figure}
\begin{centering}
\includegraphics[width=0.8\textwidth]{qq_10000sim10000_ls.eps}
\par\end{centering}

\begin{centering}
\includegraphics[width=0.8\textwidth]{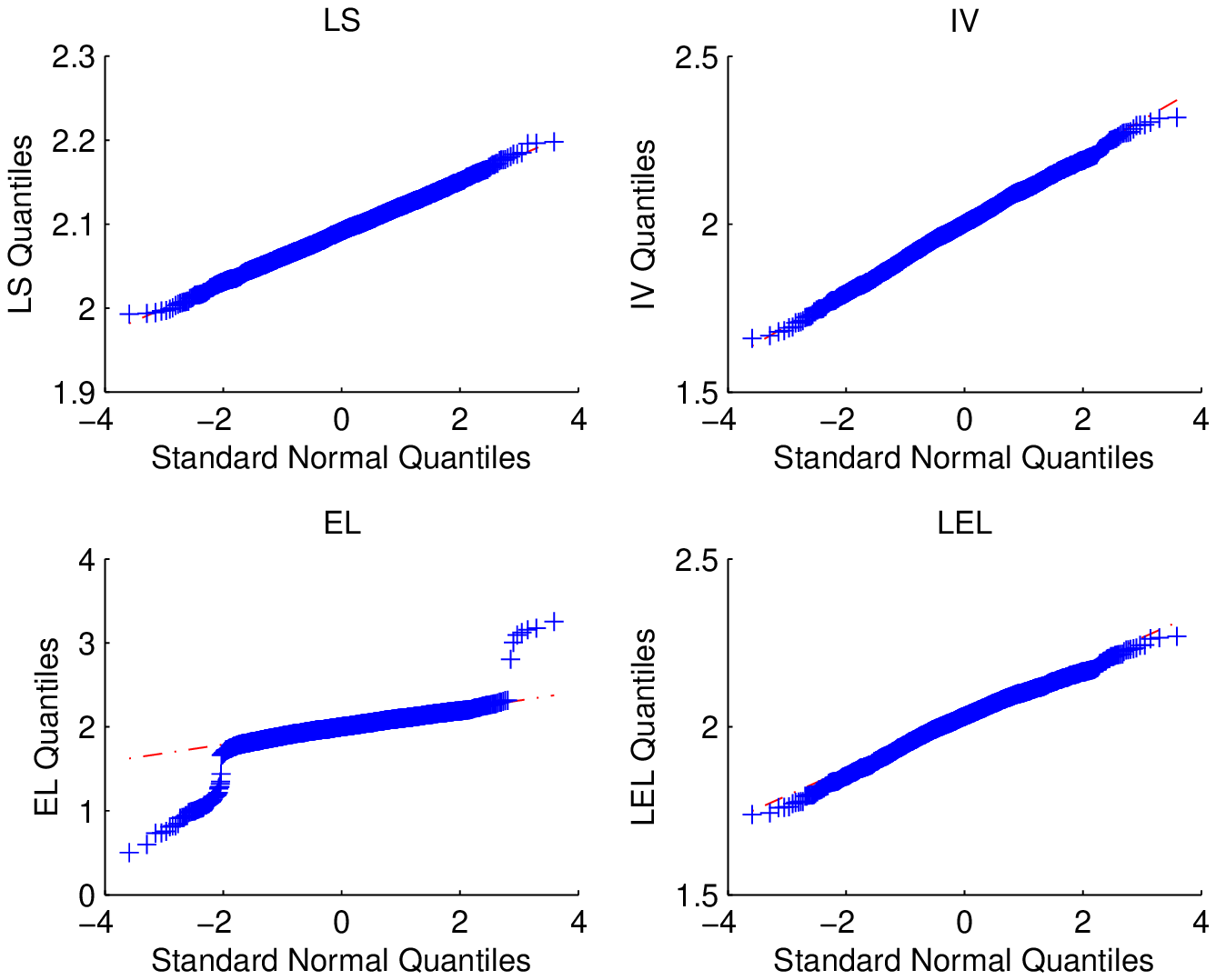}
\par\end{centering}

\caption{\label{fig:qqplot-iv}QQ plot (Densities of estimators). Up: case
(III). Down: case (IV).}
\end{figure}

\begin{figure}
\begin{centering}
\includegraphics[width=0.8\textwidth]{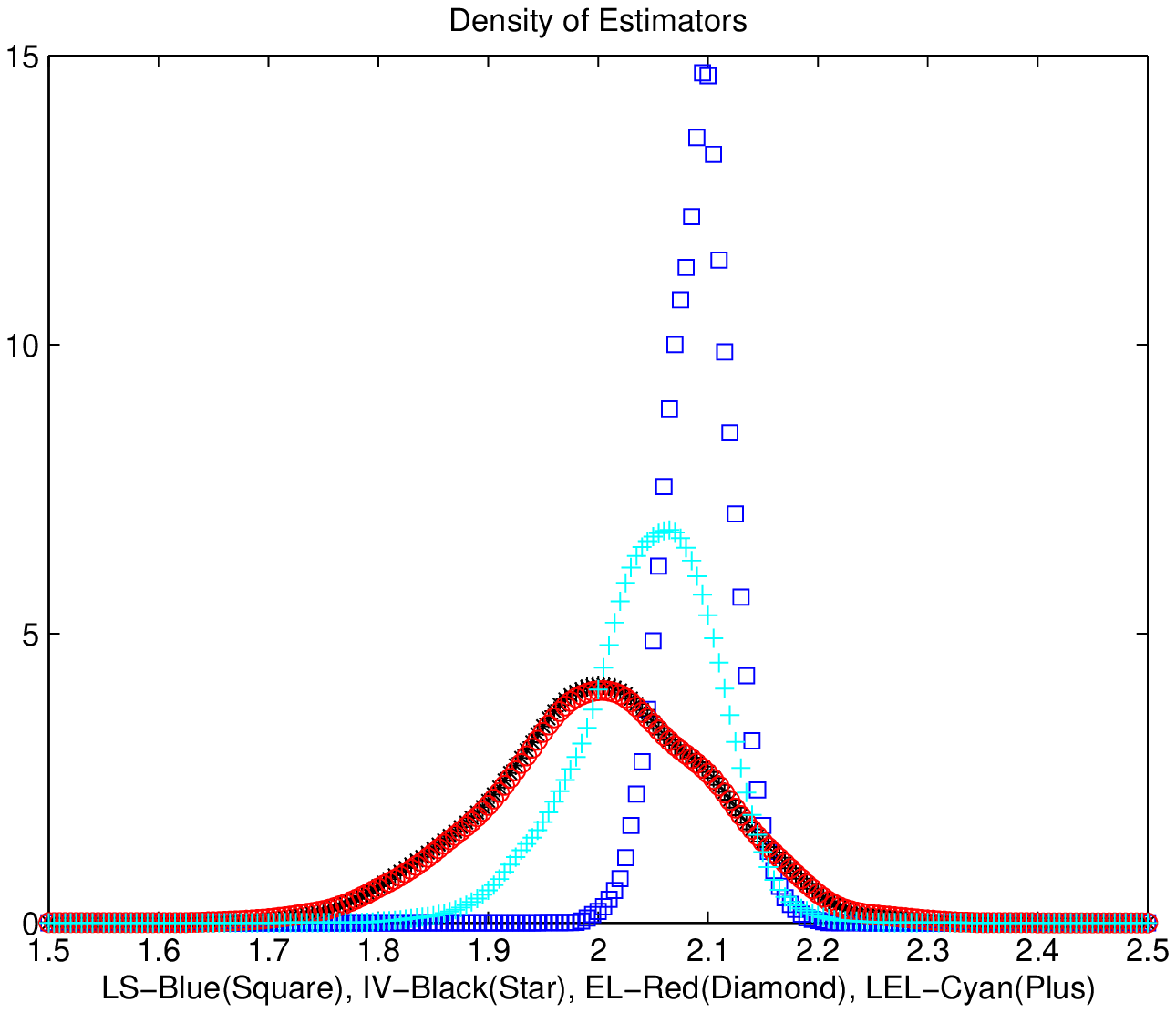}
\par\end{centering}

\begin{centering}
\includegraphics[width=0.8\textwidth]{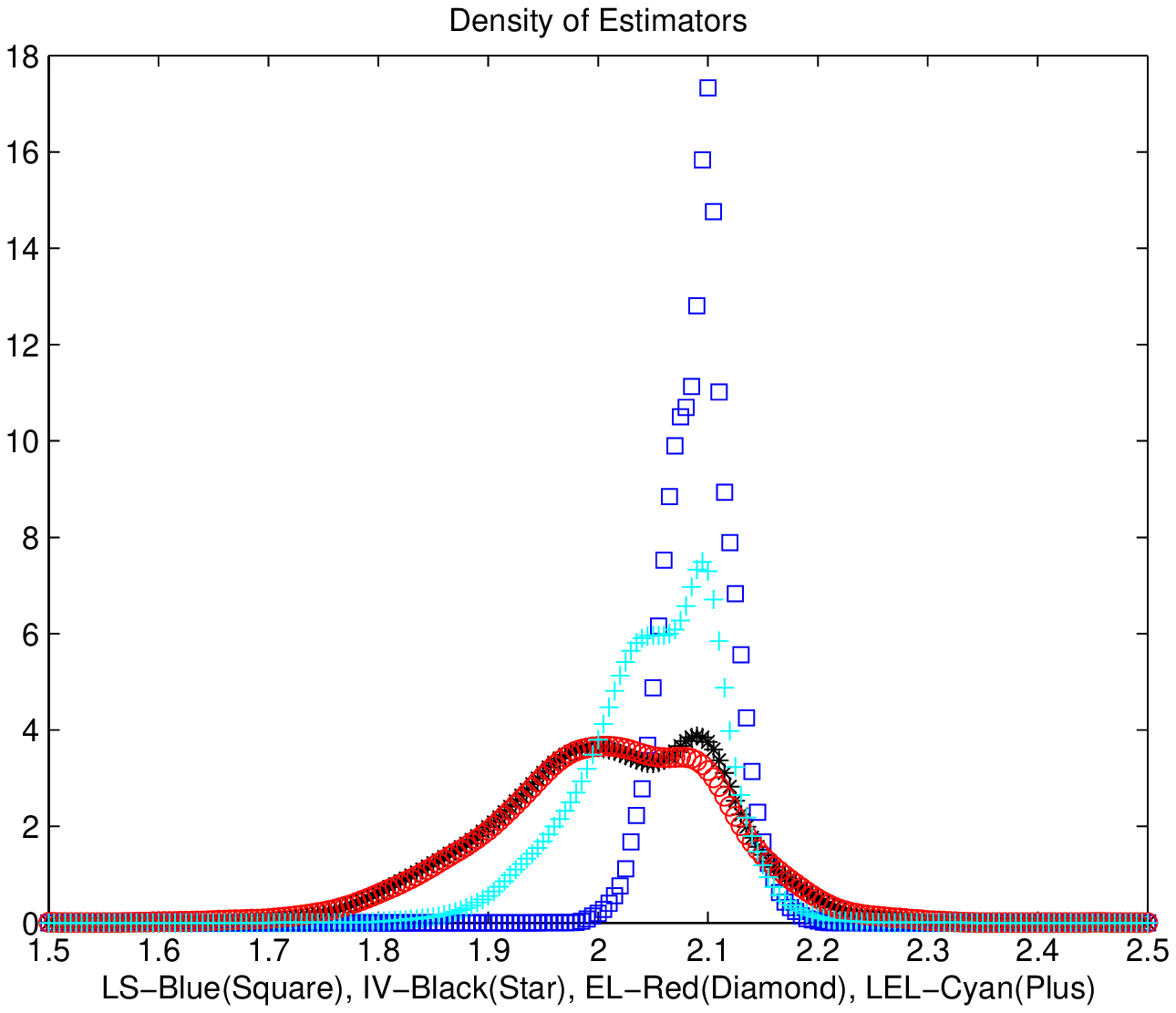}
\par\end{centering}

\caption{\label{fig:Density-plot-LS}Density plot of four estimators: LS (Blue
square), IV (Black star), EL (Red diamond) LEL (Cyan plus). Up: case
(I). Down: case (II). }
\end{figure}

\begin{figure}
\begin{centering}
\includegraphics[width=0.8\textwidth]{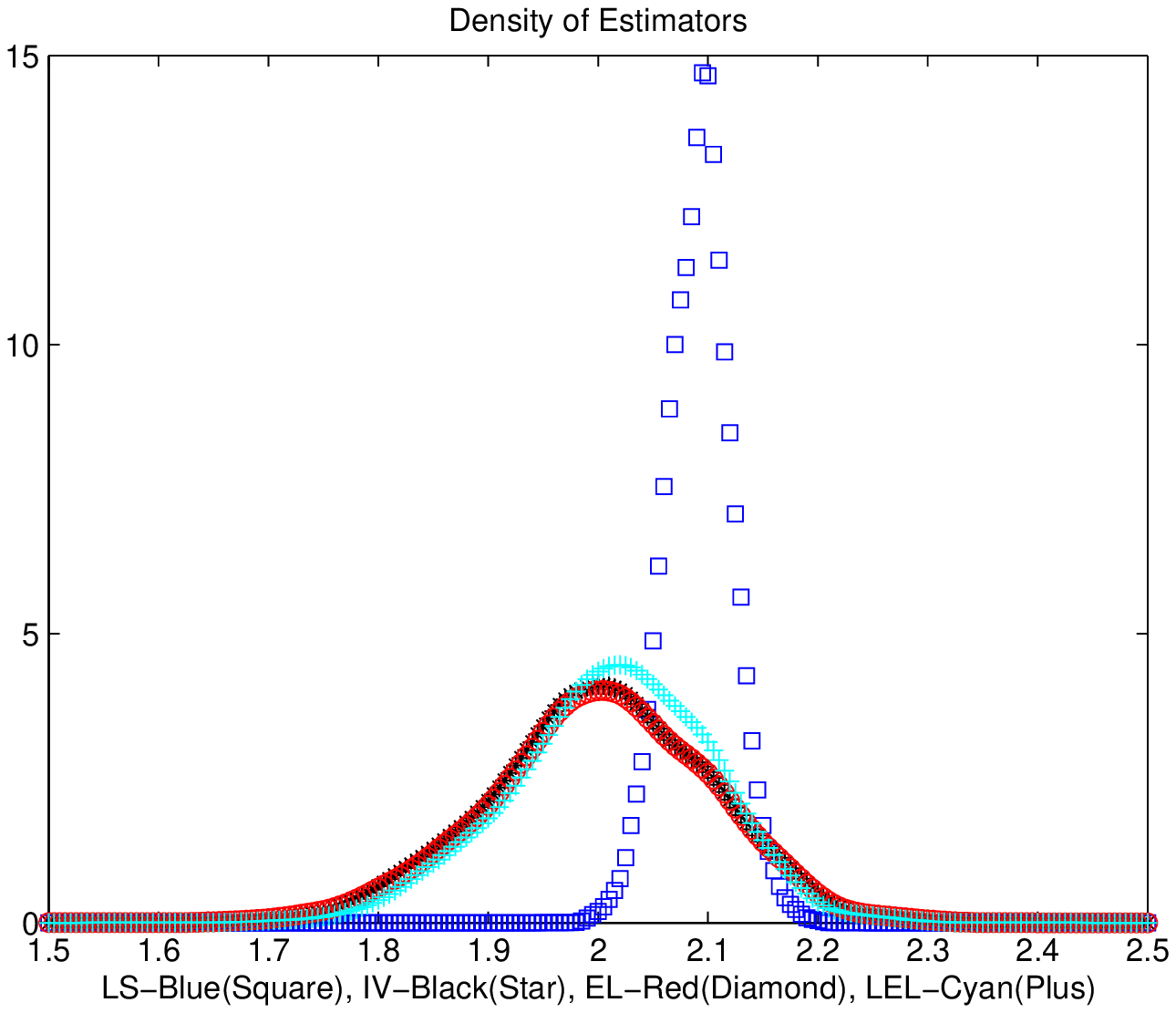}
\par\end{centering}

\begin{centering}
\includegraphics[width=0.8\textwidth]{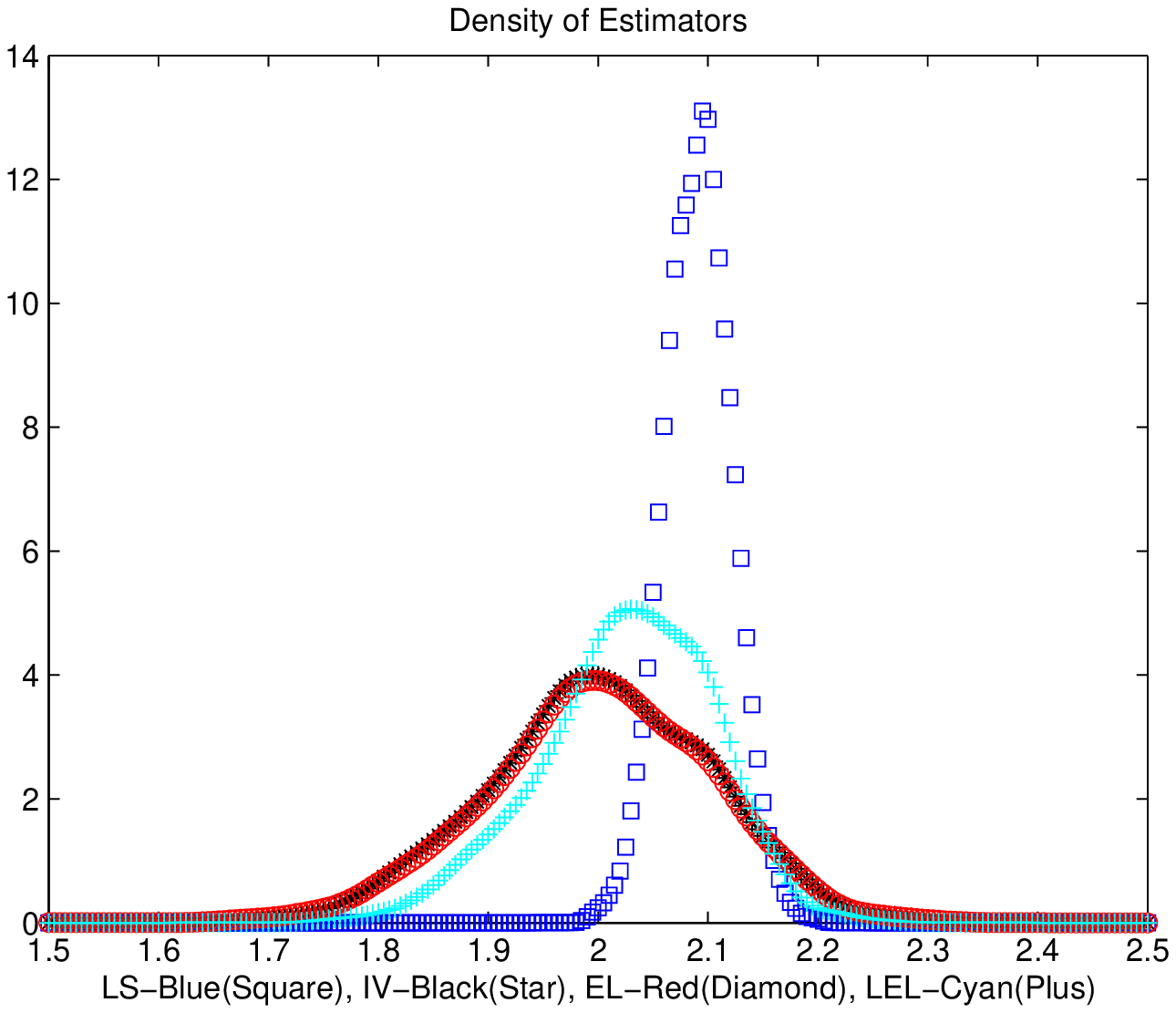}
\par\end{centering}

\caption{\label{fig:Density-plot-iv}Density plot of four estimators: LS (Blue
square), IV (Black star), EL (Red diamond) LEL (Cyan plus). Up: case
(III). Down: case (IV). }
\end{figure}

The true value of $\theta$ is $2$. A consequence of the mildly misspecified
constraint is that the moment-based estimators, IV and EL, are also
biased but not as serious as LS. The bias of LS is caused by the correlation
between $\varepsilon$ and $\mathbf{u}_{n}$. Due to the endogenous
problem, LS is always biased. The mild misspecified moment restriction
leads to the small biases in IV and EL. We will use LS or EL as the
auxiliary estimator of the local method. In Table \ref{tab:Estimation-results},
we show the estimation results for four cases: contamination percentage
$0.5\%$ ($0.005\%$) with $10$ ($10000$) degree of contamination
with LS as an auxiliary estimator; contamination level $0.01\%$ with
$10$ and $10000$ degree of contamination with IV as an auxiliary
estimator. The mean and the median of LS, IV and EL coincide with
our expectation: a large bias in LS; a relative small bias in IV and
EL. The level of bias in local method lies in-between. If one uses
LS as the auxiliary estimator, then the bias of the local method is
slightly larger than the case of using IV as the auxiliary estimator.
However, among the four estimators, local EL attains the lowest mean
square error (MSE) in all four cases. From the Q-Q plots in Figure
\ref{fig:qqplot-ls} and \ref{fig:qqplot-iv}, local EL is closer
to the normal shape than EL. The density plots in Figure \ref{fig:Density-plot-LS}
and \ref{fig:Density-plot-iv} show that the distribution of local
EL is more concentrated in case (I) and (II) but its mean location
is closer to the true value in case (III) and (IV).

\subsection{A Nonlinear Experiment}

We construct the moment restriction for a short-term interest rate
model. \citet{CLKS1992} show that the model can be nested within
the following equations:
\begin{align*}
r_{t+1}-r_{t}= & \alpha+\beta r_{t}+\varepsilon_{t+1},\\
\varepsilon_{t+1}= & \sigma r_{t}^{\gamma}u_{t},
\end{align*}
where $u_{t}$ is a normal white noise with zero mean and unit variance.
$\alpha$, $\beta$, $\gamma$, and $\sigma$ are the parameters of
the model. In this experiment, the contamination is introduced so
that the distribution of $u_{t}$ , $P_{u_{t}}$, is a mixture such
that $(1-c)P_{u_{t}^{(1)}}+cP_{u_{t}^{(2)}}.$ As in the linear case,
$u_{t}^{(1)}\sim\mathcal{N}(0,1)$, $u_{t}^{(2)}\sim\mathcal{N}(L,1)$
and $c$ denotes the contaminated percentage. to a small value so that
the model is not significantly misspecified. 

Since we have four parameters, we construct the following four moments:
\[
m_{t}(\theta)=\left[\begin{array}{c}
\varepsilon_{t+1}\\
\varepsilon_{t+1}^{2}-\sigma^{2}r_{t}^{2\gamma}
\end{array}\right]\otimes\left[\begin{array}{c}
1\\
r_{t}
\end{array}\right]=\left[\begin{array}{c}
\varepsilon_{t+1}\\
\varepsilon_{t+1}r_{t}\\
\varepsilon_{t+1}^{2}-\sigma^{2}r_{t}^{2\gamma}\\
(\varepsilon_{t+1}^{2}-\sigma^{2}r_{t}^{2\gamma})r_{t}
\end{array}\right],
\]
where $\mathbb{E}[m_{t}(\theta)]=0$. The sample moment restrictions
are 

\[
\frac{1}{T}\left[\begin{array}{c}
\sum_{t=1}^{T}((r_{t+1}-\alpha-\beta r_{t})\\
\sum_{t=1}^{T}((r_{t+1}-\alpha-\beta r_{t})r_{t}\\
\sum_{t=1}^{T}((r_{t+1}-\alpha-\beta r_{t})^{2}-\sigma^{2}r_{t}^{2\gamma}\Delta t)\\
\sum_{t=1}^{T}((r_{t+1}-\alpha-\beta r_{t})^{2}-\sigma^{2}r_{t}^{2\gamma}\Delta t)r_{t}
\end{array}\right].
\]
A consequence of the mildly misspecified constraint is that both GMM
and EL are slightly biased. The biasness is caused by the contaminated
$u_{t}$. Thus the auxiliary estimators of our local method are biased.
In this model, we will only use EL as the auxiliary estimator. 

\begin{figure}
\begin{centering}
\includegraphics[width=0.8\textwidth]{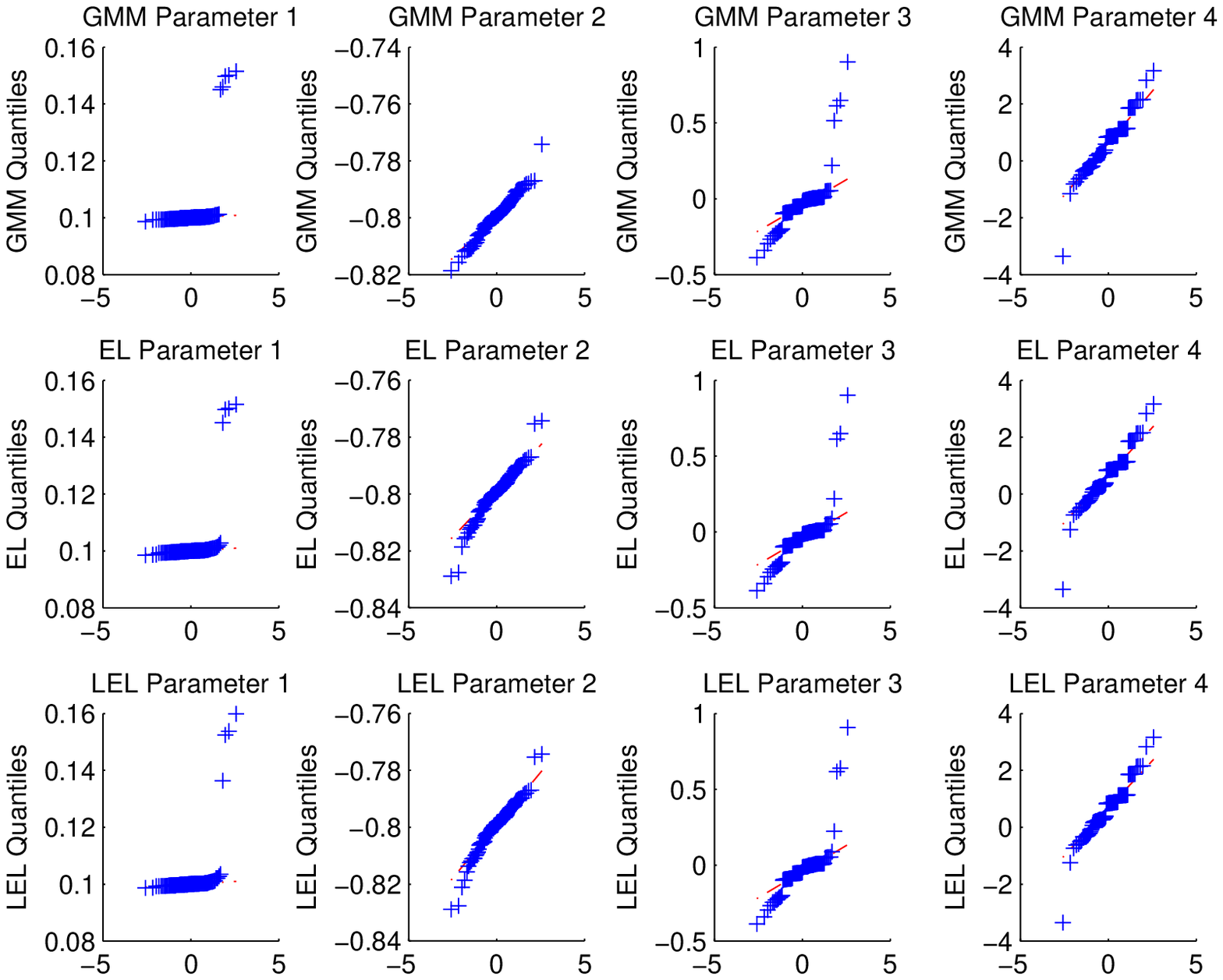}
\par\end{centering}

\caption{\label{fig:nonlinear-qqplot-1}QQ plot (Densities of estimators).
Case (I).}
\end{figure}

\begin{figure}
\begin{centering}
\includegraphics[width=0.8\textwidth]{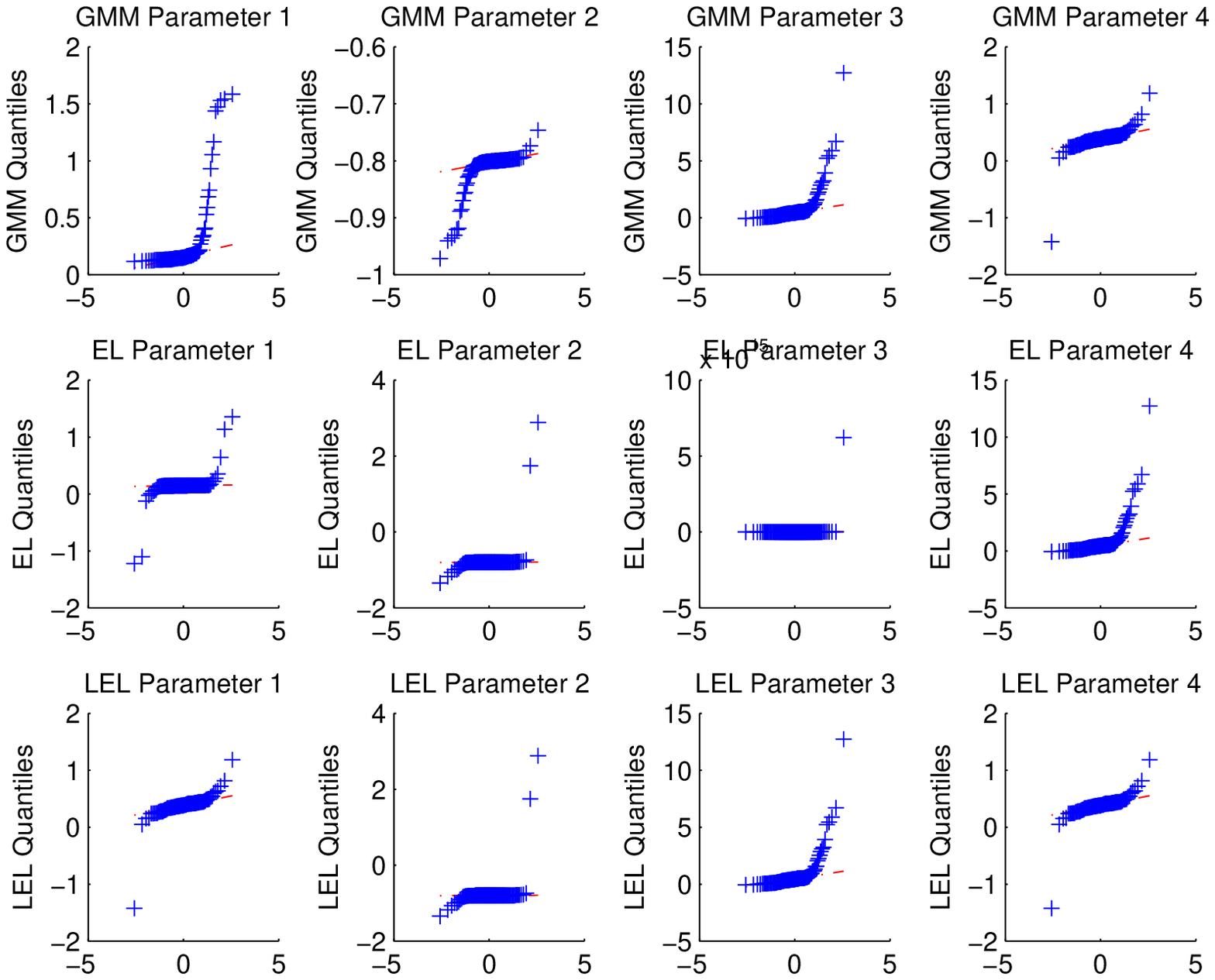}
\par\end{centering}

\caption{\label{fig:nonlinear-qqplot-2}QQ plot (Densities of estimators).
Case (II).}
\end{figure}

We restrict the contaminated level to the moderate level by setting
$L=1000$. In the experiment, we select $c$ to be $0.001\%$ and
$0.1\%$. Table \ref{tab:nonlinear-Estimation-results} shows that
the local result again lies in-between the alternative global results.
In both cases, local EL reduces root of MSE of EL. But in the small
contamination case (I), local EL is not as good as GMM because GMM
over-performs EL. While in case (II), local EL becomes a better alternative.
For estimates of each parameter, one can refer to the Q-Q plots in
Figure \ref{fig:nonlinear-qqplot-1} and \ref{fig:nonlinear-qqplot-2}.

\begin{table*}
$ $

$\:$\vspace{6pt}

\caption{\label{tab:nonlinear-Estimation-results}}

$\:$\vspace{6pt}

\centering{}{\scriptsize{}}%
\begin{tabular}{ccccccc}
\multicolumn{7}{c}{{\footnotesize{Nonlinear Model}}}\tabularnewline
\hline 
\hline 
 & \multicolumn{3}{c}{{\scriptsize{c\%=0.001\%, L=1000 (case I)}}} & \multicolumn{3}{c}{{\scriptsize{c\%=0.1\%, L=1000 (case II)}}}\tabularnewline
\hline 
{\scriptsize{Method}} & {\scriptsize{RMSE}} & {\scriptsize{IQR}} & {\scriptsize{MAD}} & {\scriptsize{RMSE}} & {\scriptsize{IQR}} & {\scriptsize{MAD}}\tabularnewline
{\scriptsize{GMM}} & {\scriptsize{0.105864}} & {\scriptsize{0.053354}} & {\scriptsize{0.053786}} & {\scriptsize{2.613145}} & {\scriptsize{0.464843}} & {\scriptsize{4.027867}}\tabularnewline
{\scriptsize{EL}} & {\scriptsize{0.106643}} & {\scriptsize{0.052433}} & {\scriptsize{0.053688}} & {\scriptsize{2.984457}} & {\scriptsize{0.509011}} & {\scriptsize{4.120982}}\tabularnewline
{\scriptsize{Local EL}} & {\scriptsize{0.106532}} & {\scriptsize{0.052420}} & {\scriptsize{0.053711}} & {\scriptsize{2.607393}} & {\scriptsize{0.467948}} & {\scriptsize{4.087653}}\tabularnewline
\hline 
\end{tabular}
\end{table*}

\section{Conclusion}

We propose a new local EL method. We discuss its construction and
derive theoretical properties. The construction is based on the infinite
divisibility property; to the best of our knowledge, this feature has not yet
been applied to EL. When the implied probability of EL is embedded
in the infinitely divisible class, the log-likelihood ratio admits
a local representation. Our local estimator is built on the basis
of this representation. The consistency, local asymptotic normality,
and asymptotic optimality of this estimator have been established.
We apply the estimate method to two simulated experiments that require
weaker regularity conditions for the estimator. The simulation results
show that the local method reduces MSE from its auxiliary estimators.

\newpage

\clearpage
\vspace{6pt}

\clearpage

\appendix

\section{Proof of Theorems}

\subsection*{Proof of Theorem \ref{thm:consistency}}

The Lagrangian of EL is
\[
L=\sum_{i=1}^{n}\log(np_{i})-n\lambda^{T}\sum_{i=1}^{n}p_{i}m_{i}(\theta)-\gamma\left(\sum_{i=1}^{n}p_{i}-1\right),
\]
 where $\lambda$ and $\gamma$ are Lagrange multipliers. Setting
the partial derivative of $L$ w.r.t $p_{i}$ equal to zero will give
$\gamma=n$ and the implied probability $\tilde{p}_{i}=1/(\gamma+n\lambda_{n}^{T}m_{i}(\theta))$.
By the implicit function theorem, the partial derivative of $\sum_{i=1}^{n}\log\tilde{p}_{i}$
w.r.t $\lambda$ gives a function $\Upsilon(\cdot,\cdot)$ of $\lambda_{n}$
and $\theta$ such that 
\begin{align}
 & \mathbf{\frac{\partial\sum\log\tilde{p}_{i}}{\partial\lambda}}:=\Upsilon(\lambda_{n},\theta)=0,\label{eq:A-1-1}\\
\Longrightarrow & \frac{1}{n}\sum_{i=1}^{n}\frac{m_{i}(\theta)}{1+\lambda_{n}^{T}m_{i}(\theta)}=\sum_{i=1}^{n}\tilde{p}_{i}(\theta)m_{i}(\theta)\notag
\end{align}
where $\lambda_{n}$ is unique for fixed $n$ and $\theta$. Note
that $\Upsilon(\lambda_{n},\theta)=0$ for $\forall\theta\in\Theta$
and $\theta$ is continuous hence $\Upsilon(\cdot)$ is continuous
in $\theta$. By the continuity of $m(X,\theta)$ and the representation
of $\Upsilon(\cdot)$, we know that $\lambda_{n}$ is also continuous
on $\theta.$ The proof of the uniqueness of $\lambda(\theta)$ is
as follows: because the set $\Gamma(\theta)=\lim_{n\rightarrow\infty}\cap_{i=1,\dots,n}\{\lambda|1+\lambda^{T}m(X_{i},\theta)>1/n\}$
is convex if it does not vanish, the function of $\log p$ is strictly
concave on $\lambda$, so $\lambda(\theta)$ exists and is unique.

With these, the properties of likelihood ratio are shown in as follows.
Equation (\ref{eq:A-1-1}) can be re-written as 
\begin{align*}
\frac{1}{n}\sum_{i=1}^{n} & \left[1-\frac{\lambda_{n}^{T}m_{i}(\theta)}{1+\lambda_{n}^{T}m_{i}(\theta)}\right]m_{i}(\theta)=0\\
\Longrightarrow\frac{1}{n}\sum_{i=1}^{n}m_{i}(\theta) & =\frac{1}{n}\sum_{i=1}^{n}\frac{m_{i}(\theta)\lambda_{n}^{T}m_{i}(\theta)}{1+\lambda_{n}^{T}m_{i}(\theta)}\\
 & =\underset{(*)}{\underbrace{\left[\sum_{i=1}^{n}\tilde{p}_{i}(\theta)m_{i}(\theta)m_{i}(\theta)^{T}\right]}\lambda_{n}}.
\end{align*}
Condition \ref{con:1-1} (v) states that $n^{-1}\sum_{i}^{n}m_{i}(\theta)m_{i}(\theta)^{T}$
is positive definite, let $\mathbf{c}$ be larger than any eigenvalue
of $n^{-1}\sum_{i}^{n}m_{i}(\theta)m_{i}(\theta)^{T}$ and let $v$
be the corresponding eigenvector. The convex combination of $m_{i}(\theta)m_{i}(\theta)^{T}$
over $\{\tilde{p}_{i}(\theta)\}$ in $(*)$ is bounded by $v^{T}\mathbf{c}v$.
Let $E_{v}=v^{T}\mathbf{c}v$. According to condition \ref{con:1-1}
(iv), $m_{i}(\theta)$ has an envelop function $b(\theta)$ such that
$\lim\inf_{\theta}|m(\theta,X)|/b(\theta)\geq1$, then 
\[
\lim_{n\rightarrow\infty}|\lambda_{n}|/b^{'}(\theta)\geq1
\]
for any $\theta$ where $b^{'}(\theta)=b(\theta)/E_{v}$.

Let's first prove the existence of $\Lambda(\theta)$: 
\begin{align}
\lim_{n\rightarrow\infty}\int\log\frac{1}{n} & \frac{n}{1+\lambda(n,\theta)^{T}m(x,\theta)}dP(x)\label{eq:A-1-2}\\
=\mathbb{E}\lim_{n\rightarrow\infty}\log\frac{1}{1+\lambda(n,\theta)^{T}m(X,\theta)} & =\mathbb{E}\log\frac{1}{1+\lambda(\theta)^{T}m(X,\theta)}=\Lambda(\theta).\notag
\end{align}
 The first convergence is by the LLN and the second equation is obtained
by the dominated convergence Theorem, since $[1+\lambda(\theta)^{T}m(x,\theta)]^{-1}$
is bounded and $\lambda(\theta)$ exists.

Next we prove the continuity of $\Lambda(\theta)$. The envelop functions
$b^{'}(\theta)$ and $b(\theta)$ are integrable and continuous (Condition
\ref{con:1-1}), $\lambda(\theta)^{T}m(X,\theta)$ is bounded by a
continuous function. Thus $\Lambda(\theta)$ is continuous and is
bounded by an envelop function $b^{''}(\theta)=\max(b^{'}(\theta),b(\theta))$
such that
\begin{equation}
\sup_{\theta}\left\Vert \Lambda_{n}(\theta)-\Lambda(\theta)\right\Vert /b^{''}(\theta)<1\label{eq:newGenCondition}
\end{equation}

Now prove the identifiability of EL estimation. Choose a compact set
$\Theta_{c}\subset\Theta$ such that for given $\epsilon$
\[
\sup_{\theta\in\Theta_{C}}|\Lambda(\theta)|/b^{''}(\theta)\geq1-\epsilon.
\]
 By (\ref{eq:newGenCondition}), LLN applied to $\Lambda_{n}(\theta)$
implies
\begin{align*}
\sup_{\theta}\frac{\left\Vert \Lambda_{n}(\theta)-\Lambda(\theta)\right\Vert }{b''(\theta)} & <\frac{\left\Vert \Lambda_{n}(\theta)\right\Vert -\left\Vert \Lambda(\theta)\right\Vert +2\left\Vert \Lambda(\theta)\right\Vert }{b^{''}(\theta)}-\epsilon\\
 & <\frac{\left\Vert \Lambda_{n}(\theta)\right\Vert -\left\Vert \Lambda(\theta)\right\Vert +2\sup_{\theta\in\Theta_{C}}|\Lambda(\theta)|}{b^{''}(\theta)}-\epsilon\\
 & <\frac{\left\Vert \Lambda_{n}(\theta)-\Lambda(\theta)\right\Vert +2\sup_{\theta\in\Theta_{C}}|\Lambda(\theta)|}{b^{''}(\theta)}<1-3\epsilon
\end{align*}
The first inequality uses triangle inequality, the second one uses
supremum property, and the third one uses triangle inequality again.
Therefore
\begin{align*}
\left|\Lambda_{n}(\theta)-\Lambda(\theta)\right| & \leq(1-3\epsilon)b^{''}(\theta)\\
\leq\frac{1-3\epsilon}{1-\epsilon}\sup_{\theta\in\Theta_{C}}|\Lambda(\theta)| & \leq(1-\delta)\sup_{\theta\in\Theta_{C}}|\Lambda(\theta)|
\end{align*}
for $\forall\theta\in\Theta_{c}$. This inequality implies
\[
\sup_{\theta\in\Theta_{c}}\left|\Lambda_{n}(\theta)\right|\leq\mathbf{\sup_{\theta\in\Theta_{c}}}|\Lambda(\theta)|+\epsilon
\]
asymptotically for any $\theta\in\Theta_{c}$. Thus if $\theta_{0}\in\Theta_{c}$,
then
\[
\{T_{n}\subset\Theta_{c}\}\subset\left\{ \sup_{\theta\in\Theta_{c}}\Lambda_{n}(\theta)\leq\Lambda(\theta_{0})+o_{p}(1)\right\} ,
\]
where the probability of the event on the right side converges to
one as $n\rightarrow\infty$. Because the compact set $\Theta$ could
be shrinking to an arbitrary neighborhood of $\theta_{0}$, the EL
estimator $T_{n}$ is consistent.

\subsection*{Proof of Theorem \ref{thm:expression}}

Before proving the theorem, we need to introduce a relation for univariate
Gaussian families. For any pair of Gaussian measures in $\mathcal{G}_{\Theta}=\{G_{\theta},\theta\in\Theta\}$,
$G_{\theta}\subset\mathcal{E}_{\theta}$, there will be an expression
to relate both of them as follows:
\begin{equation}
dG_{\theta}=\exp\left[\left\langle Y_{\vartheta},\theta\right\rangle -\frac{1}{2}\|\theta\|^{2}\right]dG_{\vartheta},\label{eq:GuassianExp}
\end{equation}
where $\vartheta,\theta\in\Theta$. The bilinear product in this expression
is $\left\langle Y_{\vartheta},\theta\right\rangle =\int_{0}^{1}Y_{\vartheta}(t)G_{\theta}(dt)$
where $Y_{\vartheta}$ is a univariate Gaussian process. This is a
random variable (functional integral or Wiener integral) with mean
zero and variance $\|\theta\|^{2}\leq\infty$ %
\footnote{This expression is called weak form expression and is often used for
generalizing Gaussian processes. %
}. If $dG_{\theta}$ and $dG_{\vartheta}$ are defined as (\ref{eq:GuassianExp}),
the integral of $(dG_{\theta}/dG_{\vartheta})^{1/2}$ w.r.t $G_{\vartheta}$
will has a linear quadratic representation.
\begin{proof}
\citet[Proposition 4.1]{LeCamYang2000} show that the affinity between
two Poissonized $d\tilde{P}_{\theta},d\tilde{P}_{\vartheta}$ is 
\[
\int\sqrt{d\tilde{P}_{\theta}d\tilde{P}_{\vartheta}}=\exp\left\{ -\frac{1}{2}\|\theta-\vartheta\|^{2}\right\} .
\]
Since \citet[Theorem 17.5]{GnedenkoKolmogorov1968} show that finite
many number of Poisson type measures can approximate any infinitely
divisible family and EL is embedded in an infinitely divisible family,
we know the above expression is applicable over here. The Hellinger
affinity for Gaussian family is 
\[
\int\sqrt{dG_{\theta}dG_{\vartheta}}=\int\exp\left[\frac{1}{2}\langle Y_{\vartheta},\vartheta+\theta\rangle-\frac{1}{4}(\|\theta\|^{2}+\|\vartheta\|^{2})\right]dG_{\vartheta}.
\]
The Gaussian property of $\langle Y_{\vartheta},\vartheta+\theta\rangle$implies
that $\exp\left[\frac{1}{2}\langle Y_{\vartheta},\vartheta+\theta\rangle\right]$
is log-normal distributed, then by log-normal property there is:
\[
\int\exp\left[\frac{1}{2}\langle Y_{\vartheta},\vartheta+\theta\rangle\right]dG_{t}=\exp\left(\frac{1}{8}\|\theta+\vartheta\|^{2}\right).
\]
Because only metric distance is going to be studied in $\int\sqrt{dG_{\theta}dG_{\vartheta}}$,
we attach a Hilbert space to $\mathcal{G}$. The parallelogram identity
for Hilbert space induces
\[
\|\theta+\vartheta\|^{2}+\|\theta-\vartheta\|^{2}=2\left(\|\theta\|^{2}+\|\vartheta\|^{2}\right),
\]
so
\[
2\left(\|\theta\|^{2}+\|\vartheta\|^{2}\right)+\|\theta+\vartheta\|^{2}=-\|\vartheta-\theta\|^{2}
\]
Therefore, $\sqrt{dG_{\theta}dG_{\vartheta}}=\exp(-\|\theta-\vartheta\|^{2}/8)$
 is isometric to $\int\sqrt{d\tilde{P}_{\vartheta}d\tilde{P}_{\theta}}=\exp(-\|\theta-\vartheta\|^{2}/2)$.
If Fubini's theorem holds, the expression 
\[
2\log\int\left(\frac{d\tilde{P}_{\theta}}{d\tilde{P}_{\vartheta}}\right)^{\frac{1}{2}}d\tilde{P}_{\vartheta}\approx8\log\int\left(\frac{dG_{\theta}}{dG_{\vartheta}}\right)^{\frac{1}{2}}dG_{\vartheta}
\]
implies 
\[
\int\left[\log\frac{d\tilde{P}_{\theta}}{d\tilde{P}_{\vartheta}}\right]d\tilde{P}_{\vartheta}=4\int\left[\log\frac{dG_{\theta}}{dG_{\vartheta}}\right]dG_{\vartheta}
\]
so that we can use the Gaussian expression (\ref{eq:GuassianExp})
for the log-likelihood ratio process. 

By Karhunen\textendash{}Loeve Theorem \citep{Kallenberg2002}, the
Gaussian process $Y_{\theta}$ can be expressed as
\[
Y_{\theta}=\sum_{j=1}^{\infty}\xi_{j}\mathbf{u}_{j}(\theta)
\]
where $\{\mathbf{u}_{j}\}$ constitutes an orthonormal basis for the
Hilbert space $\mathcal{G}$ and $\xi_{j}$ are Gaussian random variables
and stochastically independent. Now let $\mathbf{u}_{j}(\cdot)=\sum_{i}^{m}\tau_{i}\mathbf{e}_{i}(\cdot)$
where $\mathbf{e}$ is a unit basis for the local parameter space
and $\tau_{i}$ are linear coefficients for $\mathbf{e}_{i}(\cdot)$.
Let $j$ indicate the index of a basis on the Hilbert space and $i$
indicate the index of a basis on the local parameter space. Then the
inner product in the Hilbert space can be expressed using local parameter
coordinates such that $\langle Y_{\vartheta},\theta\rangle=\sum_{i}^{m}\tau_{i}\theta_{i}\langle\mathbf{e}(\vartheta),\xi\rangle=\tau^{T}(\theta\tilde{\xi})$
where $\tilde{\xi}$ is also Gaussian because of the linear property.
Let $\theta\tilde{\xi}=S_{\theta}^{'}$ and $\mathbb{E}(\theta\tilde{\xi})^{2}=K_{\theta}^{'}$,
then 
\[
\|\theta\|^{2}=\mathbb{E}[\tau^{T}(\theta\xi)]^{2}=\tau^{T}K_{\theta}^{'}\tau.
\]
From (\ref{eq:GuassianExp}), we have 
\[
\int\left[\log\frac{d\tilde{P}_{\vartheta}}{d\tilde{P}_{\theta}}\right]d\tilde{P}_{\theta}=\tau^{T}S_{\theta}-\frac{1}{2}\tau^{T}K_{\theta}\tau.
\]
where $S_{\theta}=\int S_{\theta}^{'}d\tilde{P}_{\theta}$ and $K_{\theta}=\int K_{\theta}^{'}d\tilde{P}_{\theta}$.
For a finite dimensional Gaussian vector based on $n$ realizations
Gaussian process, we have the sample counterparts $\tau_{n}$, $S_{\theta,n}$
and $K_{\theta,n}$. We conclude that the EL ratio is approximately
equal to the log-likelihood ratio of $\mathcal{G}$, which for the
sample of size $n$ is $\tau_{n}^{T}S_{\theta,n}-\tau_{n}^{T}K_{\theta,n}\tau_{n}/2$. 
\end{proof}

\subsection*{Proof of Theorem \ref{thm:localEL}}
\begin{proof}
(i) When $\theta$ is given, by equation (\ref{eq:2-1-2})
\begin{align}
\Lambda_{n}(\theta+\delta_{n}\tau_{n},\theta)= & \tau_{n}^{T}S_{\theta,n}-\frac{1}{2}\tau_{n}^{T}K_{\theta,n}\tau_{n}+o_{\tilde{p}_{\theta}}(1)\label{eq:A-4-1}\\
= & -\frac{1}{2}\left[(K_{\theta,n}^{-1}S_{\theta,n}-\tau_{n}^{T})^{T}K_{\theta,n}(K_{\theta,n}^{-1}S_{\theta,n}-\tau_{n}^{T}))\right.\notag\\
 & \left.-(S_{\theta,n}^{T}K_{\theta,n}^{-1}S_{\theta,n})\right]+o_{\tilde{p}_{\theta}}(1).\notag
\end{align}
 Similarly,
\begin{align}
\Lambda_{n}(\theta+\delta_{n}\tau_{n},\theta) & =\tau_{n}^{T}K_{n}\delta_{n}^{-1}(T_{n}-\theta_{n}^{*})-\frac{1}{2}\tau_{n}^{T}K_{n}\tau_{n}\label{eq:A-4-2}\\
 & =-\frac{1}{2}\left[(\delta_{n}(T_{n}-\theta)-\tau_{n}^{T})^{T}K_{n}(\delta_{n}(T_{n}-\theta)-\tau_{n}^{T}))\right.\label{eq:A-4-3}\\
 & \left.-(\delta_{n}(T_{n}-\theta))^{T}K_{n}(\delta_{n}(T_{n}-\theta))\right].\notag
\end{align}
The difference between (\ref{eq:A-4-1}) and (\ref{eq:A-4-2}) tends
to zero as $n\rightarrow\infty$. Non-negativity of $K_{n}$ and $K_{\theta,n}$
shows that each of the four quadratic terms in (\ref{eq:A-4-3}) and
\ref{eq:A-4-1} must be non-negative. If $S_{\theta,n}^{T}K_{\theta,n}^{-1}S_{\theta,n}$
converges to $(\delta_{n}(T_{n}-\theta))^{T}K_{n}(\delta_{n}(T_{n}-\theta))$,
then 
\begin{align*}
(\delta_{n}(T_{n}-\theta)-\tau_{n}^{T})^{T}K_{n}(\delta_{n}(T_{n}-\theta)-\tau_{n}^{T})) & \rightarrow\\
(K_{\theta,n}^{-1}S_{\theta,n}-\tau_{n}^{T})^{T}K_{\theta,n}(K_{\theta,n}^{-1}S_{\theta,n}-\tau_{n}^{T})).
\end{align*}
So one can conclude that $K_{n}\rightarrow K_{\theta,n}$ and $S_{n}\rightarrow S_{\theta,n}$.

Now consider the opposite case $(\delta_{n}(T_{n}-\theta))^{T}K_{n}(\delta_{n}(T_{n}-\theta))\nrightarrow S_{\theta,n}^{T}K_{\theta,n}^{-1}S_{\theta,n}$.
By a standard property of quadratic functions, we can have for some
positive-definite matrix $C$ 
\[
(\delta_{n}(T_{n}-\theta))^{T}K_{n}(\delta_{n}(T_{n}-\theta))+C\rightarrow S_{\theta,n}^{T}K_{\theta,n}^{-1}S_{\theta,n}
\]
Then for some vector $\Delta$ such that $\delta_{n}\Delta^{T}K_{n}\Delta\delta_{n}=C$,
there is 
\[
(\delta_{n}(T_{n}-\theta+\Delta))^{T}K_{n}(\delta_{n}(T_{n}-\theta+\Delta))\rightarrow S_{\theta,n}^{T}K_{\theta,n}^{-1}S_{\theta,n}
\]
So $T_{n}+\Delta$ is optimal estimator for $\tau_{n}$, because
\begin{align*}
(\delta_{n}(T_{n}-\theta+\Delta)-\tau_{n}^{T})^{T}K_{n}(\delta_{n}(T_{n}-\theta+\Delta)-\tau_{n}^{T})) & \rightarrow\\
(K_{\theta,n}^{-1}S_{\theta,n}-\tau_{n}^{T})^{T}K_{\theta,n}(K_{\theta,n}^{-1}S_{\theta,n}-\tau_{n}^{T})).
\end{align*}
But this contradicts with our definition of $T_{n}$.

Thus $(\delta_{n}(T_{n}-\theta))^{T}K_{n}(\delta_{n}(T_{n}-\theta))$
converges to $S_{\theta,n}^{T}K_{\theta,n}^{-1}S_{\theta,n}$. It
implies $K_{n}$ converges to $K_{\theta,n}$ in probability and $\delta_{n}(T_{n}-\theta)$
converges to $K_{\theta,n}^{-1}S_{\theta,n}$.

(ii) By Proposition \ref{pro:LAQ(iii)}, we know that clustering points
$K_{\theta}$ of $K_{\theta,n}$ are invertible. Since $\delta_{n}(T_{n}-\theta)$
converges to $K_{\theta,n}^{-1}S_{\theta,n}$, the limit of $\delta_{n}(T_{n}-\theta)$
is $K_{\theta}^{-1}S_{\theta,n}$. The Gaussian variable $S_{\theta,n}$
is second moment bounded. So the term $\delta_{n}(T_{n}-\theta)$
is bounded in probability.

(iii) We know the DQM condition implies (\ref{eq:2-1-2}), thus the
linear-quadratic equation (\ref{eq:2-1-2}) may coincide with $S_{n}$
and $K_{n}$ by (i). The log-likelihood process can be rewritten as
a centered log-likelihood process $\Xi_{n}(\cdot)$ plus a shift item
$b_{n}(\cdot)$:
\begin{align*}
\delta_{n}\Lambda_{n}(\theta,\vartheta)(x)= & \overset{\Xi_{n}(\theta)}{\overbrace{\frac{1}{n}\delta_{n}\sum_{i=1}^{n}\log\frac{\tilde{p}_{\theta}}{\tilde{p}_{\vartheta}}(x_{i})-\int\log\frac{\tilde{p}_{\theta}}{\tilde{p}_{\vartheta}}(x)dP_{0}}}\\
 & +\underset{b_{n}(\theta)}{\underbrace{\int\log\frac{\tilde{p}_{\theta}}{\tilde{p}_{\vartheta}}(x)dP_{0}}}+o_{p}(1).
\end{align*}
Let $\delta_{n}=n^{-1/2}$. Given fixed $\lambda(\cdot)$ values in
the constraint of equation (\ref{eq:2-1-0}), Theorem \ref{thm:expression}
says that $\log\frac{\tilde{p}_{\theta}}{\tilde{p}_{\vartheta}}(x_{i})$
in $\Xi_{n}(\eta)$ can be replaced by a linear quadratic formulae
w.r.t. $\tau_{n}$, namely $\log\frac{\tilde{p}_{\theta}}{\tilde{p}_{\vartheta}}(x_{i})$
belongs to a smooth functional class $\mathcal{C}^{2}$. Therefore
the process $\theta\mapsto\Xi_{n}(\theta)$ is an empirical process
and $\Xi_{n}(\theta)\rightsquigarrow\Xi(\theta)$ by Donsker's Theorem,
see \citet[Example 19.9]{Vaart1998} where $\Xi(\theta)$ is a Gaussian
process. Note that $\Xi(\theta)$ has mean $\int\Xi(\theta)dP_{0}=0$
and covariance kernel $\mathbb{E}\Xi^{2}(\theta)$ under $P_{0}$.The
log-normal property implies that $\mathbb{E}\exp[\Xi(\theta)+b(\theta)]=1$
with the expectation taken under $P_{0}$ and $b(\theta)=\lim_{n\rightarrow\infty}b_{n}(\theta)$.
Log normal property of $\exp\Xi(\cdot)$ gives $b(\theta)=-(1/2)\mathbb{E}\Xi^{2}(\theta)$.
By Proposition \ref{pro:expansion} and equation (\ref{eq:2-1-2}),
we can show that
\begin{align*}
\Xi_{n}(\theta)= & S_{\theta,n}\\
b(\theta)= & -\frac{1}{2}K_{\theta},
\end{align*}
 and when $\theta=\theta_{0}$
\begin{align*}
\Xi_{n}(\theta_{0})= & \mathbb{E}\left[\frac{\partial m(X,\theta_{0})}{\partial\theta}^{T}\left(\mathbb{E}m(X,\theta_{0})m(X,\theta_{0})^{T}\right)^{-1}\right]\delta_{n}\sum_{i}^{n}m_{i}(\theta_{0})\\
b(\theta_{0})= & -\frac{1}{2}\mathbb{E}\frac{\partial m(X,\theta_{0})}{\partial\theta}^{T}\left(\mathbb{E}m(X,\theta_{0})m(X,\theta_{0})^{T}\right)^{-1}\mathbb{E}\frac{\partial m(X,\theta_{0})}{\partial\theta}.
\end{align*}

\end{proof}

\subsection*{Proof of Theorem \ref{thm:minimax}}

Discussion: The proof follows the strategies of \citet[Proposition 8.6][]{Vaart1998}
and \citet[Theorem 6.1][]{LeCamYang1990}. The difficulty comes from
the expectation conditional on the local parameter $\tau$. Note that
the measure $\mathcal{M}$ has not yet been specified. If one can
in Bayesian fashion give a prior distribution on $\mathcal{M}$, then
what we need to study is the posterior distributions given this ``local
prior measures''. In fact, the $\delta_{n}$-sparse condition already
implies that for arbitrary priors, the corresponding posteriors concentrate
on the small shrinking neighborhood of $\theta_{0}$. 
\begin{proof}
First look at the population log-likelihood ratio
\begin{align}
\Lambda(\theta+\tau,\theta)= & -\frac{1}{2}\left[(K_{\theta}^{-1}S_{\theta}-\tau)^{T}K_{\theta}(K_{\theta}^{-1}S_{\theta}-\tau)\right.\notag\\
 & \left.-(S_{\theta}^{T}K_{\theta}^{-1}S_{\theta})\right]+o_{\tilde{p}_{\theta}}(1).\notag
\end{align}
which implies that the term $(K_{\theta}^{-1}S_{\theta}-\tau)^{T}K_{\theta}(K_{\theta}^{-1}S_{\theta}-\tau))$
is $\chi^{2}$ distributed. The quadratic form of a Gaussian variable
$\xi$, $\xi^{T}\xi$, can generate exactly the same distribution.
As Theorem \ref{thm:expression} shows that the approximation of Gaussian
family is feasible. For any value of $\theta$, there will be such
a $\xi_{\theta}$ whose distribution is equivalent to $K_{\theta}^{-1}S_{\theta}-\tau$
and has the variance $K_{\theta}^{-1/2}$. Then we have the expression
\[
\tau=K_{\theta}^{-1}S_{\theta}-\xi_{\theta},
\]
which shows that $\tau$ consists of two Gaussian variables $K_{\theta}^{-1}S_{\theta}$
and $\xi_{\theta}$. Thus we are able to impose a Gaussian structure
on the measure $\mathcal{M}$.

Now we can look at the expectation $\min(b,\mathbb{E}[W(Z_{n}-\tau)|\theta_{0}+\delta_{n}\tau])$
which is bounded by $b$. Since both ``prior'' and ``posterior''
concentrate around $\theta_{0}$ and are Gaussian, the updating information
only occurs for covariance matrix. Let $\tau$ be a Gaussian random
variable centered at $0$ with inverse covariance $\Gamma$. The conjugate
property indicates the posterior of $\tau$ can be written as: 
\[
Z_{n}=\delta_{n}^{-1}(\tilde{T}_{n}-\theta_{0})=(K_{n}+\Gamma)^{-1/2}K_{n}\delta_{n}^{-1}(T_{n}-\theta_{0}),
\]
especially when $\Gamma=0$, $Z_{n}=\delta_{n}^{-1}(T_{n}-\theta_{0})$.
By Anderson's Lemma%
\footnote{For a symmetric distribution, shifting an integral function of it
to a new position will product higher expected value, see \citet[Lemma 8.5]{Vaart1998}.%
}, for bounded $W$, there is
\[
\mathbb{E}[W(Z_{n}-\tau)|\theta_{0}+\delta_{n}\tau]\geq\mathbb{E}[W(Z_{n})|\theta_{0}+\delta_{n}\tau].
\]
Since $K_{n}\delta_{n}^{-1}(T_{n}-\theta_{0})\sim\mathcal{N}(0,I)$,
the lower bound of $\mathbb{E}[W(Z_{n}-\tau)|\theta_{0}+\delta_{n}\tau]$
is
\[
\mathbb{E}\left\{ W\left[(K_{n}+\Gamma)^{-1/2}\times\mathcal{N}(0,I)\right]|K_{n}+\Gamma\right\} .
\]
The measure of $\theta_{0}+\delta_{n}\tau$ is replaced by $K_{n}+\Gamma$
because of the Gaussian property, namely the update of covariance
matrix. Note that $K_{n}$ and $\Gamma$ are independent with $\mathcal{N}(0,I)$.
With the condition $K_{n}\rightsquigarrow K_{\theta}$ in $\tilde{P}_{\theta}$
law, the limit becomes $\mathbb{E}\left\{ W\left[(K_{\theta}+\Gamma)^{-1/2}\times\mathcal{N}(0,I)\right]\right\} .$

When $c$ is very large, the probability of normal prior $|\tau|>c$
is small enough thus
\begin{align*}
\lim\inf_{n}\sup_{|\tau|\leq c}\mathbb{E}\left\{ W\left[(K_{n}+\Gamma)^{-1/2}\times\mathcal{N}(0,I)\right]\right\}  & \geq\\
\mathbb{E}\left\{ W\left[(K_{\theta}+\Gamma)^{-1/2}\times\mathcal{N}(0,I)\right]\right\} -\Delta
\end{align*}
for small enough $\Delta.$ Especially, when $\Gamma$ go to zero
or say the measure $\mathcal{M}$ degenerates to a point eventually,
$Z_{n}=\delta_{n}^{-1}(T_{n}-\theta_{0})$ obtains the lower bound
$\mathbb{E}[W(K{}_{\theta}^{-1/2})\times\mathcal{N}(0,I)]$. If $W=1$
and $K_{\theta}=K$, by Theorem \ref{thm:localEL}(iii) we achieve
the efficient bound of semi-parametric estimators. 
\end{proof}

\section{Other Technical Details }

\subsection*{Poisson Approximation for Arbitrary Infinitely Divisible Families}

Let $\phi(t)$ and $\phi_{n}(t)$ be the characteristic functions
of distributions in $\mathcal{E}$ and $\mathcal{E}_{n}$. By the
infinitely divisible property, $\phi(t)=[\phi_{n}(t)]^{n}$ or $\phi_{n}(t)=[\phi(t)]^{1/n}$.
Two characteristic functions have the following relation: 
\begin{align*}
n(\phi_{n}(t)-1)=n(\sqrt[n]{\phi(t)}-1) & =n\left(e^{\frac{1}{n}\log\phi(t)}-1\right)\\
=n\left(1+\frac{1}{n}\log\phi(t)+o(\frac{1}{n})-1\right) & \rightarrow\log\phi(t),
\end{align*}
or say $\exp(n(\phi_{n}(t)-1))\rightarrow\phi(t)$. The concrete construction
of characteristic function in $\mathcal{E}_{\theta,n}$ depends on
the discrete Fourier transform of $\Lambda(X,\theta)$ on $j$ segments
e.g. $\inf\Lambda(X)<c_{1}<c_{2}<\dots<c_{j}<\sup\Lambda(X)$ which
implies that 
\[
\lim_{j\rightarrow\infty}\sum_{k=1}^{j}a_{k}(i)e^{itc_{k}}=\int e^{it\Lambda(X)}dF_{n}=\phi_{n}(t),
\]
where $a_{n}(k)=n(F_{n}(c_{k})-F_{n}(c_{k-1}))$ is the Fourier coefficient%
\footnote{The Stieltjes sum, a discrete version of stochastic integral.%
} and $F_{n}$ is the measure for $\Lambda_{n}(\theta)$. Combined
with the expression above, one can see that a characteristic function
of finite many number of Poisson measures (compound Poisson measures)
approximates $\phi(t)$:
\begin{equation}
\exp\sum_{i=1}^{j}(na_{i})\left(e^{it\Lambda(x_{i},\theta)}-1\right)\rightarrow\phi(t)\label{eq:poissonCharacter}
\end{equation}
where $j\rightarrow\infty$ and $\{na_{i}\}_{i=1,\dots,j}$ converges
to a measure. To see the argument of (\ref{eq:poissonCharacter}),
let $V(\cdot)$ be a Poisson process (a random measure) with Poisson
parameter $\gamma$ such that $\mathbb{E}V(\mathcal{A})=\gamma(\mathcal{A})$
for a set $\mathcal{A}$. For any function $v$ in infinite divisible
family, the characteristic function of $v$ is $\phi(t)=\exp\{\int(e^{itv}-1)d\gamma\}$. 

The approximation can be viewed as constructing a new family which
approximately equals the infinite divisible $\mathcal{E}_{\theta}$.
Firstly select a Poisson variable $\nu$ (again a random measure)
such that $\mathbb{E}\nu(\Lambda(X))=1$ for any log-likelihood ratio
$\Lambda(X)$ and then carry out $n$-draws from the direct product
$\otimes_{i=1,\dots,\nu}\mathcal{E}_{\theta,i}$, $\nu$ copies $\mathcal{E}_{\theta,i}$.
The result is called a poissonized family.

\subsection*{Derivation of Equation (\ref{eq:step2-2})}

Since $\mathbb{E}[\exp(\log(dG_{i}/d\mu))]=1$, then we have 
\[
\mathbb{E}\exp\left[L(i)+\mathbb{E}\log(\frac{dG_{i}}{d\mu})\right]=\left[\mathbb{E}e^{L(i)}\right]\cdot e^{\mathbb{E}\log(\frac{dG_{i}}{d\mu})}=1
\]
By the log-normal property, $\mathbb{E}\exp L(i)=e^{\frac{1}{2}K(i,i)}$,
we have 
\[
e^{\frac{1}{2}K(i,i)}\cdot e^{\mathbb{E}\log(\frac{dG_{i}}{d\mu})}=1\Longleftrightarrow\mathbb{E}\left[\log(\frac{dG_{i}}{d\mu})\right]=-\frac{1}{2}K(i,i)
\]
thus we have (\ref{eq:step2-1}). For $\mathbb{E}\exp[L(\theta)+L(\vartheta)]$,
we have $2K(\theta,\vartheta)$. Combining $2K(\theta,\vartheta)$
and $K(i,i)$ gives us (\ref{eq:step2-2}).

\subsection*{Proof of Proposition \ref{pro:expansion}}

The proof is based on Taylor expansions. Note that 
\begin{equation}
m(x,\theta_{0}+\delta_{n}\tau)=m(x,\theta_{0})+\delta_{n}\frac{\partial m(x,\theta_{0})}{\partial\theta^{T}}\tau+o_{p}(\delta_{n}^{2}).\label{eq:A-2-1}
\end{equation}
 Let $\theta\in\{\theta||\theta-\theta_{0}|\leq|\tau|\delta_{n}\}$,
$|\tau|$ is a vector with elements equal to their absolute values.
The result 
\[
\lambda_{n}(\theta)=\left(\sum_{i=1}^{n}[m_{i}(\theta)m_{i}(\theta)^{T}]/n\right)^{-1}\sum_{i=1}^{n}m_{i}(\theta)/n+o_{p}(n^{-1/2})
\]
holds uniformly for $\theta$ in a neighborhood of $\theta_{0}$,
see the proofs in \citet[Lemma 1]{QinLawless1994} or \citet[Theorem 2.2]{Owen2001}.
For the empirical log-likelihood at $\theta$, by noting that $\lambda_{n}^{T}m_{i}$
is close to zero and using a second order approximation for $\log(1+\lambda_{n}^{T}m_{i}),$we
obtain: \begin{align*}
\sum_{i=1}^{n}\log\tilde{p}_{\theta}= & \sum_{i=1}^{n}\left[\lambda_{n}(\theta)^{T}m_{i}(\theta)-\frac{1}{2}\left(\lambda_{n}(\theta)^{T}m_{i}(\theta)m_{i}(\theta)^{T}\lambda_{n}(\theta)\right)\right]\\
 & -n\log n+o_{p}(1).
\end{align*}
The remainder term is based on bounding $\sum_{i=1}^{n}\left(\lambda_{n}^{T}m_{i}\right)^{3}$
for which \citet{Owen1990} showed in Lemma 3 that it is of order
$o_{p}(1)$. Note that his $\gamma_{i}$ is our $\lambda_{n}^{T}m_{i}(\theta)$.
Note that 
\[
\lambda_{n}(\theta)^{T}m_{i}(\theta)=\left(\sum_{i=1}^{n}\frac{m_{i}(\theta)}{n}\right)^{T}\left[\sum_{i=1}^{n}\frac{1}{n}\left(m_{i}(\theta)m_{i}(\theta)^{T}\right)\right]^{-1}m_{i}(\theta)
\]
 and after summation equals the squared term: 
\begin{align*}
\sum_{i=1}^{n} & \lambda(\theta)^{T}m_{i}(\theta)m_{i}(\theta)^{T}\lambda_{n}(\theta)=\\
\left(\sum_{i=1}^{n}\frac{m_{i}(\theta)}{n}\right)^{T} & \left[\sum_{i=1}^{n}\frac{1}{n}\left(m_{i}(\theta)m_{i}(\theta)^{T}\right)\right]^{-1}\left(\sum_{i=1}^{n}\frac{m_{i}(\theta)}{n}\right).
\end{align*}
 So adding these two terms we obtain: 
\begin{align*}
\sum_{i=1}^{n}\log\tilde{p}_{\theta} & =\frac{1}{2}\left(\sum_{i=1}^{n}\frac{m_{i}(\theta)}{n}\right)^{T}\left[\sum_{i=1}^{n}\frac{1}{n}\left(m_{i}(\theta)m_{i}(\theta)^{T}\right)\right]^{-1}\\
 & \times\left(\sum_{i=1}^{n}\frac{m_{i}(\theta)}{n}\right)-n\log n+o_{p}(1).
\end{align*}
It implies: 
\begin{align*}
2\sum_{i=1}^{n}\log\frac{\tilde{p}_{\theta_{0}+\delta_{n}\tau}}{\tilde{p}_{\theta_{0}}}(x_{i})=\left(\frac{1}{n}\sum_{i=1}^{n}m_{i}(\theta_{0}+\delta_{n}\tau)\right)^{T} & \times\\
\left(\frac{1}{n}\sum_{i=1}^{n}[m_{i}(\theta_{0}+\delta_{n}\tau)m_{i}(\theta_{0}+\delta_{n}\tau)^{T}]\right)^{-1} & \sum_{i=1}^{n}m_{i}(\theta_{0}+\delta_{n}\tau)-\\
\left(\frac{1}{n}\sum_{i=1}^{n}m_{i}(\theta_{0})\right)^{T}\left(\frac{1}{n}\sum_{i=1}^{n}[m_{i}(\theta_{0})m_{i}(\theta_{0})^{T}]\right)^{-1} & \sum_{i=1}^{n}m_{i}(\theta_{0})+o_{p}(1).
\end{align*}
It follows from the approximation of $\lambda$ above. Using equation
(\ref{eq:A-2-1}) we can further simplify the terms involving $\theta+\delta_{n}\tau$.
We obtain for the middle term:
\begin{multline*}
\frac{1}{n}\sum_{i=1}^{n}[m_{i}(\theta_{0}+\delta_{n}\tau)m_{i}(\theta_{0}+\delta_{n}\tau)^{T}]=\frac{1}{n}\sum_{i=1}^{n}\left[[m_{i}(\theta_{0})m_{i}(\theta_{0})^{T}]+\right.\\
\delta_{n}\tau\left(\frac{\partial m_{i}(\theta_{0})}{\partial\theta^{T}}\right)^{T}m_{i}(\theta_{0})+\frac{(\delta_{n}\tau)^{2}}{4}\left(\frac{\partial m_{i}(\theta_{0})}{\partial\theta^{T}}\right)^{T}\frac{\partial m_{i}(\theta_{0})}{\partial\theta^{T}}\left.+o_{p}(\delta_{n}^{3})\right]\\
=\frac{1}{n}\sum_{i=1}^{n}[m_{i}(\theta_{0})m_{i}(\theta_{0})^{T}]+\frac{1}{n}\delta_{n}O_{p}(n^{1/2})+o_{p}(\delta_{n}^{2})+o_{p}(\delta_{n}^{3}).
\end{multline*}
With the big bracket becoming
\begin{align*}
 & n\left[\frac{1}{n}\sum_{i=1}^{n}m_{i}(\theta_{0})+\frac{1}{n}\sum_{i=1}^{n}\delta_{n}\frac{\partial m_{i}(\theta_{0})}{\partial\theta^{T}}\tau\right]^{T}\left(\frac{1}{n}\sum_{i=1}^{n}[m_{i}(\theta_{0})m_{i}(\theta_{0})^{T}]\right)^{-1}\\
 & \times\left[\frac{1}{n}\sum_{i=1}^{n}m_{i}(\theta_{0})+\frac{1}{n}\sum_{i=1}^{n}\delta_{n}\frac{\partial m_{i}(\theta_{0})}{\partial\theta^{T}}\tau\right]\\
=n & \left[\frac{1}{n}\sum_{i=1}^{n}m_{i}(\theta_{0})+\delta_{n}\mathbb{E}\frac{\partial m_{i}(\theta_{0})}{\partial\theta^{T}}\tau+\delta_{n}O(n^{-1/2}(\log\log n)^{1/2})\right]^{T}\\
 & \times\left(\mathbb{E}\left(m(x,\theta_{0})m(x,\theta_{0})^{T}\right)\right)^{-1}\\
 & \times\left[\frac{1}{n}\sum_{i=1}^{n}m_{i}(\theta_{0})+\delta_{n}\mathbb{E}\frac{\partial m_{i}(\theta_{0})}{\partial\theta^{T}}\tau+\delta_{n}O(n^{-1/2}(\log\log n)^{1/2})\right]\\
=2 & \delta_{n}\mathbb{E}\frac{\partial m(x,\theta_{0})}{\partial\theta^{T}}\tau\left(\mathbb{E}\left(m(x,\theta_{0})m(x,\theta_{0})^{T}\right)\right)^{-1}\frac{1}{n}\sum_{i=1}^{n}m_{i}(\theta_{0})+\\
 & \delta_{n}^{2}\mathbb{E}\frac{\partial m(x,\theta_{0})}{\partial\theta^{T}}\tau\left(\mathbb{E}\left(m(x,\theta_{0})m(x,\theta_{0})^{T}\right)\right)^{-1}\mathbb{E}\frac{\partial m(x,\theta_{0})}{\partial\theta}\tau+\\
 & \frac{1}{n}\sum_{i=1}^{n}m_{i}(\theta_{0})\left(\mathbb{E}\left(m(x,\theta_{0})m(x,\theta_{0})^{T}\right)\right)^{-1}\frac{1}{n}\sum_{i=1}^{n}m_{i}(\theta_{0})+o_{p}(\delta_{n}^{3})
\end{align*}
 where $O(n^{-1/2}(\log\log n)^{1/2})$ is used to bound the difference
of the sample average and the expectation of a random vector. Thus
the local EL is
\[
2\sum_{i=1}^{n}\log\frac{\tilde{p}_{\theta_{0}+\delta_{n}\tau_{n}}}{\tilde{p}_{\theta_{0}}}(x_{i})=\delta_{n}\tau_{n}^{T}A_{1}++\frac{1}{2}\delta_{n}^{2}\tau_{n}^{T}A_{2}\tau_{n}^{T}+o_{p}(1)
\]
where
\begin{align*}
A_{1}= & \mathbb{E}\frac{\partial m(X,\theta_{0})}{\partial\theta}^{T}\left(\mathbb{E}m(X,\theta_{0})m(X,\theta_{0})^{T}\right)^{-1}\sum_{i=1}^{n}m_{i}(\theta_{0}),\\
A_{2}= & \mathbb{E}\frac{\partial m(X,\theta_{0})}{\partial\theta}^{T}\left(\mathbb{E}m(X,\theta_{0})m(X,\theta_{0})^{T}\right)^{-1}\mathbb{E}\frac{\partial m(X,\theta_{0})}{\partial\theta^{T}}.
\end{align*}
Note that $O(n^{-1/2}(\log\log n)^{1/2})\times\delta_{n}\sum_{i=1}^{n}m_{i}(\theta_{0})=o_{p}(1)$
and 
\[
\lim_{n\rightarrow\infty}A_{n}\cdot\sum_{i=1}^{n}[m_{i}(\theta_{0}+\delta_{n}\tau)-m_{i}(\theta_{0})]/n=o_{p}(1)
\]
with $A_{n}=\sum_{i=1}^{n}m_{i}(\theta_{0})\left(\mathbb{E}m(x,\theta_{0})m(x,\theta_{0})^{T}\right)^{-1}$
by the continuity of $m_{i}(\theta)$.

\subsection*{Proof of Proposition \ref{pro:LAQ(iii)}}

To prove $K_{\theta}$ is invertible, we will prove $K_{\theta}$
is almost surely positive definite. Le Cam's first Lemma implies that
\begin{equation}
\mathbb{E}\exp\left[\tau^{T}S_{\theta}-\frac{1}{2}\tau^{T}K_{\theta}\tau\right]=1.\label{eq:A-3-1}
\end{equation}
Because (\ref{eq:A-3-1}) holds for all $\tau$, we can use a symmetrized
method to simplify (\ref{eq:A-3-1}). For a given value $\tau$ and
$-\tau$, we have 
\[
\mathbb{E}\left\{ \exp\left[\tau^{T}S_{\theta}-\frac{1}{2}\tau^{T}K_{\theta}\tau\right]+\exp\left[-\tau^{T}S_{\theta}-\frac{1}{2}\tau^{T}K_{\theta}\tau\right]\right\} =2.
\]
 By $\cosh\tau^{T}S_{\theta}=(\exp\tau^{T}S_{\theta}+\exp(-\tau^{T}S_{\theta}))/2$,
we have
\begin{equation}
\mathbb{E}[(\cosh\tau^{T}S_{\theta})\exp(-\tau^{T}K_{\theta}\tau/2)]=1.\label{eq:A-20}
\end{equation}
Assume there is some $\tau$ such that $\tau^{T}K_{\theta}\tau$ is
negative, then 
\begin{align}
\mathbb{E}\left[\mathbb{I}_{\{\tau^{T}K_{\theta}\tau>0\}}(\cosh\tau^{T}S_{\theta})\exp(-\tau^{T}K_{\theta}\tau/2)\right]\label{eq:A-19}\\
\leq\mathbb{E}\left[(\cosh\tau^{T}S_{\theta})\exp(-\tau^{T}K_{\theta}\tau/2)\right] & =1\notag
\end{align}
 where $\mathbb{I}_{\{\cdot\}}$ is an indicator function. However,
since
\[
\exp(-\tau^{T}K_{\theta}\tau/2)>1
\]
when $\tau^{T}K_{\theta}\tau$ is negative and $(\cosh\tau^{T}S_{\theta})>1$,
\begin{align*}
 & \underset{>0}{\underbrace{\mathbb{E}\left[\mathbb{I}_{\{\tau^{T}K_{\theta}\tau>0\}}(\cosh\tau^{T}S_{\theta})\exp(-\tau^{T}K_{\theta}\tau/2)\right]}+}\\
 & \underset{\geq1}{\underbrace{\mathbb{E}\left[\mathbb{I}_{\{\tau^{T}K_{\theta}\tau\leq0\}}(\cosh\tau^{T}S_{\theta})\exp(-\tau^{T}K_{\theta}\tau/2)\right]}}\\
 & =\mathbb{E}\left[(\cosh\tau^{T}S_{\theta})\exp(-\tau^{T}K_{\theta}\tau/2)\right]>1
\end{align*}
we have a contradiction with equation (\ref{eq:A-20}) unless the
set $\{\tau^{T}K_{\theta}\tau\leq0\}$ is empty. Therefore, $K_{\theta}$
is positive definite and hence invertible.

\section{Implementation of the Local EL in Section \ref{sec:Local-Estimation}}

The evaluation of the LEL estimator requires evaluation of $S_{n}$
and $K_{n}$. It appears reasonable to use any numerical first and
the second derivative of $\Lambda_{n}(\theta_{n}^{*}+\delta_{n}\tau,\theta_{n}^{*})$.
The matrix $u_{i}^{T}K_{n}u_{j}=\{K_{n,i,j}\}$ in Section \ref{sec:Local-Estimation}
\begin{align*}
K_{n,i,j}= & -\left\{ \Lambda_{n}[\theta_{n}^{*}+\delta_{n}(u_{i}+u_{j}),\theta_{n}^{*}]\right.\\
 & \left.-\Lambda_{n}[\theta_{n}^{*}+\delta_{n}u_{i},\theta_{n}^{*}]-\Lambda_{n}[\theta_{n}^{*}+\delta_{n}u_{j},\theta_{n}^{*}]\right\} 
\end{align*}
is a particular form of a numerical derivatives. If $u\in\mathbb{R}$,
the above expression of $K_{n,i,j}$ can be simplified to
\[
K_{n}=\frac{-\Lambda_{n}[\theta_{n}^{*}+\delta_{n}2u,\theta_{n}^{*}]+2\Lambda_{n}[\theta_{n}^{*}+\delta_{n}u,\theta_{n}^{*}]}{u^{2}}.
\]
For a fixed value of $\delta_{n}$, if we let $f(\delta_{n}u)=\Lambda_{n}[\theta_{n}^{*}+\delta_{n}2u,\theta_{n}^{*}]$,
then there is
\[
\delta_{n}^{-2}K_{n}=-\frac{[f(2\delta_{n}u)-f(\delta_{n}u)]/\delta_{n}u-[f(\delta_{n}u)-f(0)]/\delta_{n}u}{\delta_{n}u}
\]
which is a simple one-sided numerical second derivative of $f(\delta_{n}u)$
at $u=0$, multiplied by $-1$. Note that 
\[
\lim_{\delta_{n}u\rightarrow0}\frac{\Lambda_{n}[\theta_{n}^{*}+\delta_{n}u,\theta_{n}^{*}]-0}{\delta_{n}u}
\]
define the derivative of $\Lambda_{n}[\theta_{n}^{*}+\delta_{n}u,\theta_{n}^{*}]$
at $\theta_{n}^{*}$. In our implementation, instead of using the
expression $(f(\delta_{n}u)-f(0))/\delta_{n}u$, we will focus on
the derivative form of $\Lambda_{n}[\theta_{n}^{*}+\delta_{n}u,\theta_{n}^{*}]$.
While $\lambda_{n}(\theta_{n}^{*})$ in $\tilde{p}_{\theta_{n}^{*}}$
cannot be attained as a closed form expression, we will use the Romberg
method to handle this difficulty.

The whole implementation of LEL follows the definition in Section
\ref{sec:Local-Estimation}.

Step 1. Find an auxiliary estimate $\theta_{n}^{*}$ using LS or IV.

Step 2. can be written as an expression of a 2nd order finite difference
\begin{align*}
K_{n,i,j}=- & \left\{ \log\frac{\tilde{p}_{\theta_{n}^{*}+2\delta_{n}u}}{\tilde{p}_{\theta_{n}^{*}}}-\log\frac{\tilde{p}_{\theta_{n}^{*}+\delta_{n}u}}{\tilde{p}_{\theta_{n}^{*}}}-\log\frac{\tilde{p}_{\theta_{n}^{*}+\delta_{n}u}}{\tilde{p}_{\theta_{n}^{*}}}\right\} \\
=- & \left\{ \frac{1}{2}\cdot2\left(\log\tilde{p}_{\theta_{n}^{*}+2\delta_{n}u}-\log\tilde{p}_{\theta_{n}^{*}}\right)-\right.\\
 & \left.2\left(\log\tilde{p}_{\theta_{n}^{*}+\delta_{n}u}-\log\tilde{p}_{\theta_{n}^{*}}\right)\right\} \\
=- & [f(2\delta_{n}u)-f(\delta_{n}u)]-[f(\delta_{n}u)-f(0)]
\end{align*}
Then $(f(\delta_{n}u)-f(0))/\delta_{n}u$ can be expressed as a directional
derivative $\frac{\partial}{\partial\vec{u}}\log\tilde{p}$ evaluated
at $\theta_{n}^{*}$:
\begin{align*}
\frac{1}{\delta_{n}u}\left(\log\tilde{p}_{\theta_{n}^{*}+\delta_{n}u_{i}}-\log\tilde{p}_{\theta_{n}^{*}}\right) & \rightarrow\frac{\partial}{\partial\vec{u}}\log\tilde{p}
\end{align*}
as $\delta_{n}u\rightarrow0$. Similar argument holds for $(f(2\delta_{n}u)-f(\delta_{n}u))/\delta_{n}u$
which is the directional derivative evaluated at $\theta_{n}^{*}+\delta_{n}u$. 

Hence, the Hessian is constructed by the directional derivative. We
need to obtain the numerical value of the directional derivative $\frac{\partial}{\partial\vec{u}}\log\tilde{p}$.
Using the chain rule, a directional derivative $\frac{\partial}{\partial\theta}\log\tilde{p}_{\theta}$
can be expressed as $\partial_{\theta}(\lambda m)\times(\tilde{p}_{\theta})^{-1}$
where $\partial_{\theta}(\lambda m)=\partial_{\theta}(\lambda_{n}(\theta)m(X,\theta))$
is a numerical derivative using the Romberg method%
\footnote{Because there is no closed form expression for $\lambda$, there is
no way of obtaining analytical expression of $\partial_{\theta}(\lambda_{n}(\theta)m(X,\theta))$. %
}, see e.g. \citet{Kiusalaas2009}. 

The next task is to find a proper direction $u$. Because the direction
$u$ can be arbitrarily chosen%
\footnote{The direction $u_{i}$ and $u_{j}$ are unknown. The directional derivative
$\frac{\partial}{\partial\vec{u}}(\cdot)$ depends on $u_{i}$ and
\uline{$u_{j}$}.%
}. We simply search the direction $u$ using bisection method. 

The bisection method concerns on $\tilde{\theta}=\theta_{n}^{*}+\delta_{n}u$
such that 
\[
\lambda(\theta)\sum_{i=1}^{n}m_{i}(\theta)=-\lambda_{n}(\theta)\sum_{i=1}^{n}m_{i}(\tilde{\theta}),
\]
where $\sum_{i}m_{i}(\theta)=Z^{T}(Y-X\theta)$ in our experiment.
So the simplified expression of $\tilde{\theta}$ is
\[
\lambda(\theta)Z^{T}(Y-X\theta)=-\lambda(\theta)Z^{T}(Y-X\tilde{\theta})
\]
or $X\tilde{\theta}=(2Z^{T}Y-Z^{T}X\theta).$ Then the directional
derivative $\frac{\partial}{\partial\vec{u}}\log\tilde{p}$ can be
set to $\partial_{\tilde{\theta}}(\lambda m)(\tilde{p}_{\tilde{\theta}})^{-1}$.

The Hessian used in the implementation is
\begin{align*}
\delta_{n}^{-2}K_{n,i,j}\approx- & \left[\frac{\partial}{\partial\vec{u}}\left(\frac{\partial}{\partial\vec{u}_{1}}\log\tilde{p}+\frac{\partial}{\partial\vec{u}_{2}}\log\tilde{p}\right)\right]\\
=\partial_{\tilde{\theta}}(\lambda m) & \partial_{\theta}(\lambda m)\left[\frac{1}{(\tilde{p}_{\theta_{n}^{*}})^{2}}+\frac{1}{(\tilde{p}_{\tilde{\theta}})^{2}}\right]\partial_{\theta}(\lambda m)\partial_{\tilde{\theta}}(\lambda m).
\end{align*}

Step 3. After some rearrangement of $f(u)$, the linear term $S_{n}$
can be expressed as:
\[
\delta_{n}^{-1}S_{n}=\frac{3}{2}\frac{f(\delta_{n}u)-f(0)}{\delta_{n}u}-\frac{1}{2}\frac{f(2\delta_{n}u)-f(\delta_{n}u)}{\delta_{n}u},
\]
which is a weighted average of numerical first derivative of $f(\tau)$
at $\tau=0$ and $\tau=u$. We simply use $\partial_{\theta_{n}^{*}}(\lambda m)(\tilde{p}_{\theta_{n}^{*}})^{-1}$
to express $\delta_{n}^{-1}S_{n}$.

Step 4. Construct the adjusted estimator:
\[
T_{n}=\theta_{n}^{*}+\delta_{n}K_{n}^{-1}S_{n}=\theta_{n}^{*}+(\delta_{n}^{2}K_{n}^{-1})\times(\delta_{n}^{-1}S_{n})
\]

\bibliographystyle{authoryear}
\bibliography{\jobname}

\end{document}